\documentclass[12pt]{article}

\usepackage[bottom=0.8in,top=0.8in]{geometry}
\usepackage{amssymb,latexsym,amsmath,amsfonts}
\usepackage[mathscr]{eucal}
\usepackage{mathrsfs}
\usepackage{amsthm}
\usepackage{amsxtra}
\usepackage{txfonts}
\usepackage{float}
\usepackage[usenames,dvipsnames]{xcolor}
\usepackage{hyperref}
\hypersetup{colorlinks=true, citecolor=MidnightBlue, linkcolor=BrickRed}
\usepackage{setspace}
\usepackage{subfiles}
\usepackage{preamble2}

\newcommand{\Addresses}{{
  \bigskip
  \footnotesize

  \textsc{Department of Mathematics at Rutgers University, }
   \texttt{ogrodnikbphd@gmail.com}
}}

\title{On the Local-Global Conjecture for Commutator Traces}

\begin{document}
\author{Brooke Logan Ogrodnik\footnote{\Addresses}}

\maketitle

\abstract{

We study the trace set of the commutator subgroup of $\Gamma(2),$ a type of Local-Global problem about thin groups. We determine the local obstructions and then use the correspondence between binary quadratic forms and hyperbolic matrices to find some global obstructions. We then develop a probabilistic argument for the existence of sufficiently large admissible traces by modeling the elements as non-backtracking random walks in a 2-dimensional lattice via their homology class and word length. Finally, we investigate the number of commutators needed to represent matrices in the commutator subgroup. This is done using an algorithm of Goldstein and Turner along with utilizing properties of level $k$-Markoff type surfaces. We conjecture that any trace in the commutator subgroup of $\Gamma(2)$ can be represented with either 1 or 2 commutators. 
}

\section{Introduction}\label{ch:Intro}

The questions addressed in this paper are motivated by a conjecture of McMullen \cite{mcmullen} and a special case of this conjecture is as follows.

\begin{conj}\label{mc}
The set $\{[\overline{a_1, \ldots, a_k}]\in \mathbb{Q}(\sqrt{5}) \mid a_i\leq 2\}$ has exponential growth as $k\rightarrow \infty.$
\end{conj}

In 2018, Bourgain and Kontorovich reformulated this conjecture into a Local-Global question about thin (semi) groups \cite{BK4}. Before we show the reformulation, we give the general definition of admissibility and critical exponents.

\begin{defn}\label{admiss}
Let $\Gamma$ be a sub-semigroup of $\SL_2(\mathbb{Z})$ and $F$ a polynomial with respect to the entries of the matrices such that $F(\Gamma)\subset \mathbb{Z}$. An integer $n$ is defined as \textbf{admissible} if for all $q\geq1$, $n\in F(\Gamma) \bmod{q}$.
\end{defn}

\begin{defn}
Let $\Gamma$ be a sub-semigroup of $\SL_2(\mathbb{Z}).$ Then its \textbf{critical exponent} is defined to be the abscissa of convergence of its Poincar\'e series, $\sum_{\gamma \in \Gamma} \norm{\gamma}^{-2s}.$
\end{defn}

Using the observation that $[\overline{a_1, \ldots, a_k}]$ is fixed by $\begin{pmatrix} a_1 & 1\\1 & 0 \end{pmatrix}\cdots \begin{pmatrix} a_k & 1\\1 & 0 \end{pmatrix}$, the conjecture by McMullen can be reformulated to that of a Local-Global problem in the following way. Let $A$ be a subset of the integers. Define
\begin{equation*}
\Gamma_A:=\left<\begin{pmatrix} a & 1\\1 & 0 \end{pmatrix}\mid a\in A \right>^+ \cap \SL_2(\mathbb{Z}),
\end{equation*} (where $<g_i>^+$ is the semigroup generated by the elements $g_i$) choose $F=\Tr$, and look at the trace set $\mathscr{T}_A:=\{\Tr(\gamma) \mid \gamma\in \Gamma_A\}$ ($\mathscr{T}_A=F(\Gamma)$). If $1,2\in A$ then it is easy to see that every $n$ is admissible.

\begin{conj}[Local-Global with Multiplicity \cite{BK4}]\label{multi}
 Let $A$ be an alphabet for which the critical exponent,
$\delta_A$, exceeds 1/2. Then the set $\mathscr{T}_A$ of traces contains every sufficiently large admissible integer. Moreover,
the multiplicity for admissible $t\in [X,2X]$ is at least $$\#\{\gamma \in \Gamma_A \mid \Tr(M)=t\text{ and } \left\lVert \gamma \right\rVert<X\}>X^{2\delta_A-1-o(1)}.$$
\end{conj}

Conjecture \ref{multi} would then imply McMullen's Conjecture \ref{mc} (see \cite{BK4}). Despite this reformulation of the problem, relatively little progress has been made. In this paper we consider a related problem in a particular group that has more structure. The group is defined as follows.

 Take $\Gamma(2):=\{\gamma \in \SL_2(\mathbb{Z}) \mid \gamma \equiv I \bmod{2}\}$, the level-2 principal congruence subgroup of the modular group $\SL_2(\mathbb{Z})$, and let $\Gamma(2)'$ be its commutator subgroup. Our main objective, analogous to McMullen's conjecture, is to study the trace set of $\Gamma(2)'.$ This is interesting since the trace is related to the length of closed geodesics in the trivial homology. First, we determine the admissible values.
 
\begin{thm}[Admissibility Theorem]
The set of admissible values for the trace set of $\Gamma(2)'$ are all $t$ such that $t \bmod{256} \in \{2,18,66,146\}$ and $t \bmod{9} \in \{0,2,3,6,7\}.$
\end{thm}
\begin{conj}
For all sufficiently large admissible $t$, there exists $\gamma\in \Gamma(2)'$ such that $\Tr(\gamma)=t.$
\end{conj}

This conjecture asks for less than the what Conjecture \ref{multi} would ask in the analogous setting since we are only asking if matrices with trace $t$ exist, and not about multiplicity. We are not able to fully solve this and so in these pages we will prove some theorems about the structure of the group and give both empirical and probabilistic arguments toward the conjecture.

It is not even known that a positive proportion of the traces arise. One's first instinct may be to use the relationships between level structures of the Markoff surface and 1-commutators, elements of the form $[X,Y]$ with $X,Y\in \Gamma(2)$ (see Section \ref{ch:coms}), as studied by Ghosh and Sarnak \cite{ghosh2017integral}. That is, perhaps the traces of 1-commutators already produce a positive proportion. However, this direction seems just as, if not more, difficult since one would need to exhibit matrices in $\Gamma(2)$ with very special properties. This question about 1-commutators then opens an auxiliary direction of inquiry where one asks, for $\gamma$ in the commutator subgroup of $\Gamma(2)',$ what is the commutator width of $\gamma$? I.e. the minimal number of commutators needed to represent $\gamma.$ There are cases of different groups where such a question is undecidable \cite{notcomputable}, however Wicks showed that it is possible to compute the width in free groups \cite{MR142610}. We use the algorithm of Goldstein and Turner \cite{GoldsteinRichard1979Aotg} to do so effectively and study the distribution of widths which lead to the formulation of Conjecture \ref{comsconjecture}.

Another question one might ask is, are there admissible values that do not arise as traces? We explain a constructive algorithm (Section \ref{ch:existence}) that will tell us if a given value is or is not a trace of a matrix in the commutator subgroup of $\Gamma(2)$ modeling Gauss's reduction theory of Binary Quadratic Forms. In doing this, we identify what the failures are in this context and it seems reasonable to believe that these are all of them. 

\begin{lem}\label{lem6}
The only admissible values, $t$, such that $|t|<10^6$ and $t$ is not the trace of an element in $\Gamma(2)'$ is when $t\in \{-1006, -1726, -2558, -4718\}.$
\end{lem}

\begin{conj}
The only admissible values, $t$, such that $t$ is not the trace of an element in $\Gamma(2)'$ is when $t\in \{-1006, -1726, -2558, -4718\}.$
\end{conj}

Note that while lengths of closed geodesics correspond to $\pm \Tr(g),$ the positive values of these local-global failures are not admissible since in our group if $t$ is admissible, $-t$ is not.

In Section \ref{ch:walks} we model representatives of the $\SL_2(\mathbb{Z})$-conjugacy classes of $H$ as random walks on a 2 dimensional lattice and use this to derive a heuristic for the number of conjugacy classes in the commutator subgroup with trace $t$ using the assumption that on average these walks are of length $\log|t|.$

This research began with experimentation. All the physical code for this project can be found at \href{https://github.com/BrookeOgrodnik/CommutatorSubgroups}{https://github.com/BrookeOgrodnik/CommutatorSubgroups} and for examples and illustrations of some of the objects looked at throughout the paper, one can go to \href{https://classnumbers-and-walks.herokuapp.com/}{https://classnumbers-and-walks.herokuapp.com/} to view an interactive app.

\subsection{Notations, Observations, and Theorems}

 Trace is a function on conjugacy classes, so for this local-global problem we want to know for each admissible value if a conjugacy class exists with that trace. Let us define the following congruence subgroup $L:=\{\gamma\in \SL_2(\mathbb{Z}) \mid \gamma \bmod{8}\in \{I,5I\}\}.$ 
\begin{itemize}
  \item \textbf{Observation 1: } We note that $\Gamma(2)'\triangleleft L$ (see Lemma \ref{lem1}) and so looking at $L$ (instead of $\Gamma(2)$) reduces the search space for elements in the trivial homology.
\end{itemize}

For convenience we will also let 
\begin{equation*}
C:=\begin{pmatrix}1 & 1 \\ 0 & 1\end{pmatrix}, D:=\begin{pmatrix}1 & 0 \\ 1 & 1\end{pmatrix}, \mathcal{A}:=C^2, \text{ and } \mathcal{B}:=D^2. 
\end{equation*}
Recall that the group generated by $C$, $D$, and $-I$ is $\SL_2(\mathbb{Z})$. From here we can define the free group generated by two elements as $H:=\left<\mathcal{A},\mathcal{B}\right>$ which has the property that $\Gamma(2)=\left<H,-I\right>$. This gives us that $\Gamma(2)'= H'$ so we can (and will) talk about the two interchangeably. Then we also define $\Delta:=\{\mathcal{A}^{m_1}\mathcal{B}^{n_1}\cdots \mathcal{A}^{m_k}\mathcal{B}^{n_k} \mid k\geq 0, m_i,n_i\neq 0\}$ and note that $H=\cup_{m,n\in \mathbb{Z}} B^n \Delta A^m.$ 

\begin{itemize}
  \item \textbf{Observation 2: }Every conjugacy class of $\Gamma(2)'$ can be represented by an element in $\Delta.$
  \item \textbf{Observation 3: } In determining which traces arise in $\Gamma(2)'\triangleleft L$ we may group elements under conjugation with respect to $\SL_2(\mathbb{Z})$ since $L\triangleleft \SL_2(\mathbb{Z}).$
\end{itemize}

We will be using different versions of ``class number" throughout these pages, two of which are defined below. 
\begin{defn}\label{classnum}
Define the \textbf{$L$ class number}, $h(t)$, as the number of conjugacy classes in $\SL_2(\mathbb{Z})$ of hyperbolic matrices in $L$ with trace $t.$
\end{defn}

\begin{defn}
Define the \textbf{ $\Gamma(2)'$ class number}, $h'(t)$, as the number of conjugacy classes in $\SL_2(\mathbb{Z})$ of hyperbolic matrices in $\Gamma(2)'$ with trace $t.$
\end{defn}
We now give a lower bound (hinted at in Section \ref{algoExistence} and explained in detail in \cite{dissertation}) on the number of conjugacy classes with trace $t$ in $L$.
\begin{cor}\label{corexist}
For admissible $t,$ the $L$ class numbers satisfy $h(t)\gg_\epsilon |t|^{1-\epsilon}$ where the implied constant is ineffective.
\end{cor}
Next, we give a lower bound for the word length of $\gamma\in \Delta$ (see Definition \ref{narrowlength}).

\begin{thm}
If $g\in \Delta$ and $\gamma\neq I$, then the word length of $g$ is bounded below by $\log(|\Tr(\gamma)|/2)$.
\end{thm} This lower bound inspires the following conjecture when one models the average length of $\gamma$ as $O(\log|\Tr(\gamma)|)$.

\begin{conj}
For admissible $t,$ the $\Gamma(2)'$ class number satisfies $h'(t)\asymp h(t)/\log(|t|).$
\end{conj}
\begin{figure}[h]
\center
\includegraphics[width=60mm]{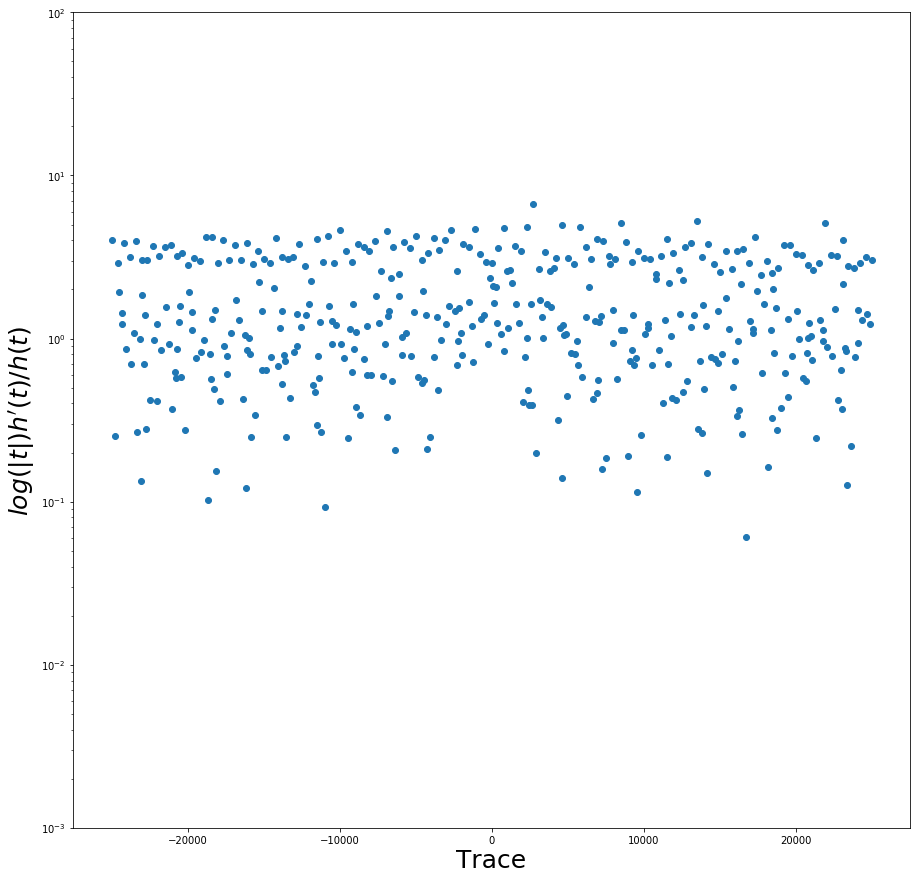}
\caption{A log plot of $(t,h'(t)\log(|t|)/h(t))$ for admissible $t$ with $|t|<25000$}
\label{ht}
\end{figure}

Finally, when investigating the width question (Section \ref{ch:coms}) we see that it is not true that every trace has a representative that is a 1-commutator. Thus, we propose this final conjecture.

\begin{conj}\label{comsconjecture}
For every $t$ in $\Tr(\Gamma(2)')$, there exists a $\gamma\in \Gamma(2)'$ such that $\Tr(\gamma)=t$ and $\gamma$ is either of commutator width 1 or 2.
\end{conj}

\section{Preliminaries}\label{ch:Pres}

In the previous section we defined $H=<\mathcal{A}, \mathcal{B}>$ and its commutator subgroup, $H'$. The following lemma contains statements that were either directly or indirectly mentioned already and the proofs are elementary and thus omitted.

\begin{lem}\label{appendixlemma}\label{obvi}The following hold:
\begin{enumerate}
    \item $H$ and $H'$ are free groups
    \item  $H/H'\cong \mathbb{Z}^2$
    \item  $H'$ is an infinitely generated group
    \item $H'=\left<[\mathcal{A}^m,\mathcal{B}^n] \mid n,m \in \mathbb{Z}\right>.$
    \item For all $g\in H,$ there exists $m,n\in \mathbb{Z}$ and $h\in H'$ such that $g=\mathcal{A}^m\mathcal{B}^nh$
    \item  $\Gamma(2)=\pm H$ and $H'= \Gamma(2)'.$
\end{enumerate}
\end{lem} 

Take $g=\mathcal{A}^{m_1}\mathcal{B}^{n_1}\cdots \mathcal{A}^{m_k}\mathcal{B}^{n_k}\in H$ where $n_i, m_i\in \mathbb{Z}$. Then a \textbf{walk} will be defined as the following set of points:
\begin{equation*}
\left\{(0,0), (m_1, 0), (m_1, n_1), (m_1+m_2, n_1), \ldots, \left(\sum m_i, \sum n_i\right)\right\}.
\end{equation*} Then let $\pi(g):=(\sum m_i, \sum n_i)$ be the homology class that $g$ falls in.

\begin{defn}\label{narrowlength}
Take $g\in H$ with $g=\mathcal{A}^{m_1}\mathcal{B}^{n_1}\cdots \mathcal{A}^{m_k}\mathcal{B}^{n_k}$ where $m_i, n_i\in\mathbb{Z}$ and $m_1$ and $n_k$ are the only values allowed to be equal to 0. Then the \textbf{word length of $g$} is $\ell(g)=\sum |m_i|+|n_i|.$ We further define the \textbf{narrow length} for $g \in \Delta$ to be $\ell_N(g)=k.$
\end{defn}

\begin{defn}
For any element $g\in \Gamma(2)'$, $g$ can be expressed as a product of $k$-commutators, $g=[g_1, h_1]\cdots [g_k, h_k]$ for some $k\in \mathbb{N}.$ The least such $k$ is called the \textbf{commutator width} of $g$. 
\end{defn}

 Clearly, $\Gamma(2)'\triangleleft \GL_2(\mathbb{Z})$ since $\Gamma(2) \triangleleft \GL_2(\mathbb{Z}).$ Thus for any $\Gamma$ such that $\Gamma(2)'\leq \Gamma \leq \GL_2(\mathbb{Z})$, it is also true that $\Gamma(2)'\triangleleft \Gamma.$


\begin{lem}\label{lem1} $\Gamma(2)'\triangleleft L\triangleleft \GL_2(\mathbb{Z})$ and furthermore, $g\in L$ implies $\Tr(g)\equiv 2 \bmod{16}.$
\end{lem}

\begin{proof}
 We can easily calculate the group $\Gamma(2)' \bmod{8}$. We know that $\Gamma(2)'$ is generated by $[\mathcal{A}^m,\mathcal{B}^n]$ and since we are looking modulo 8 we just need to check what these generators are when $0\leq n,m<3 $ and we find that it is $\{I, 5I\}$. Thus $\Gamma(2)'\leq L$ and a normal subgroup by the previous statement. One can show that $g \in \Gamma(4)$ implies that $\Tr(g)\equiv 2\bmod{16}$ and, since $L\subset \Gamma(4)$, the same holds in $L.$ Furthermore $L\triangleleft \GL_2(\mathbb{Z})$ given that $L$ is a congruence subgroup.
\end{proof} This lemma gives us the first insight to the local obstructions of the trace set, i.e. that any trace in the commutator subgroup must, at the very least, be congruent to 2 modulo 16.

\subsection{Local-Global}

A result of Epstein \cite{eps} (see also Sharp \cite{sharp}), which extended the work of Phillips and Sarnak \cite{phisar},  gives us the following (in crude form) for $\Gamma(2)'$:
\begin{thm}\label{growth}
The critical exponent for $\Gamma(2)'$ is $\delta=1.$ Equivalently, $$\#\{\gamma\in \Gamma(2)' \mid \left\lVert \gamma\right\rVert<X\}=X^{2+o(1)}$$ where $\left\lVert \cdot \right\rVert$ is any Archimedian norm.
\end{thm}

Given that $\Gamma(2)'$ is the commutator subgroup of $H$, we know it is an infinite index subgroup of $\Gamma(2)$ (by Lemma \ref{appendixlemma}) and thus an infinite index subgroup of $\SL_2(\mathbb{Z})$ and is also Zariski dense (Corollary \ref{zariskicoro}), so it is \textbf{thin}. 

\begin{defn}
An integer $n$ is \textbf{admissible} if for all $q\geq1$, $n\in \Tr(\Gamma(2)') \bmod{q}$.
\end{defn}

We will now state another version of the Local-Global conjecture as well as a lemma demonstrating the bounds that we get for free due to the structure of the group. Similar to Conjecture \ref{multi}, we define the multiplicity as $
  \mult_X(t)=\{\gamma\in \Gamma(2)'\mid \left\lVert \gamma\right\rVert<X \textit{ and has trace } t\}.$

\begin{conj}
If $t\asymp X$ is admissible then, $\mult_X(t)>X^{1-o(1)}.$
\end{conj}
Almost for free, one gets just shy of positive proportion of traces as seen in the following lemma.
\begin{lem}
\begin{equation*}
 \#\{|t|<X\mid \exists \gamma \in \Gamma(2)', \left\lVert \gamma\right\rVert<X, \Tr(\gamma)=t\}>X^{1-o(1)}.\end{equation*}
\end{lem}
\begin{proof}
First we note that $\mult_X(t)\leq X^{1+o(1)}.$ This comes from the fact that in order to construct $\begin{pmatrix}a & b\\c & d\end{pmatrix}$, there are $\ll X$ choices for $a$, and once we choose $a$, the value $d=t-a$ is determined. Finally there are $X^{o(1)}$ possibilities for $b$ and $c$ as divisors of $ad-1$.

Given Theorem \ref{growth}, we know that 
\begin{equation*}
  \sum_{|t|<X} \mult_X(t)=X^{2+o(1)}
\end{equation*} and combined with the multiplicity bound above, we use Cauchy-Schwarz to see that 
 \begin{align*}
  X^{2+o(1)}&=\sum_{|t|<X} \mult_X(t)\cdot 1_{\mult_X(t)\neq 0}
\end{align*}
And so the number of distinct traces less than $X$ is greater than $X^{1-o(1)}$.
\end{proof}

\section{The Local Theory}\label{ch:Admi}

Throughout this section and the next, $\SL_2(k):=\SL_2(\mathbb{Z}/k\mathbb{Z})$ and with this in mind, these pages are dedicated to proving the two theorems below. This section relies on work in the GitHub repository. \footnote{\href{https://github.com/BrookeOgrodnik/CommutatorSubgroups/blob/master/Mathematica/groupmodk.nb}{https://github.com/BrookeOgrodnik/CommutatorSubgroups/blob/master/Mathematica/ \newline groupmodk.nb}}
\begin{thm}\label{big} For $k$ a prime power or coprime to 6,
\begin{equation*}
  \Gamma(2)'\bmod{k}= \begin{cases} 
\{I \} &\text{for } k=1, 2,4\\
\{I, 5I \} &\text{for } k=8\\
g_1 &\text{for } k=16\\
 \{\gamma \in \SL_2(2^l) : \gamma \in \Gamma(2)' \bmod{16}\}
 &\text{for } k=2^l \text{, } l\geq 5\\
g_2 &\text{for } k=3\\
 \{\gamma \in \SL_2(3^l) : \gamma \in \Gamma(2)' \bmod{3}\} & \text{for } k=3^l \text{, } l\geq 2\\
\SL_2(k) & \text{else }
\end{cases}
\end{equation*} where 
\begin{equation*}
  g_1=\left\{I,9I,\begin{pmatrix}5 & 8 \\ 8 &13 \end{pmatrix},\begin{pmatrix}13 & 8 \\ 8 &5 \end{pmatrix}  \right\}
\end{equation*} and 
\begin{equation*}
  g_2=\left\{I,2I,\begin{pmatrix}0 & 2 \\ 1 &0 \end{pmatrix},\begin{pmatrix}1 & 1 \\ 1 &2 \end{pmatrix} ,\begin{pmatrix}1 & 2 \\ 2 &2 \end{pmatrix} ,\begin{pmatrix}0 & 1 \\ 2 &0 \end{pmatrix},\begin{pmatrix}2 &2 \\ 2 &1 \end{pmatrix}, \begin{pmatrix}2 & 1 \\ 1 &1 \end{pmatrix}\right\}.
\end{equation*}
\end{thm} 
These prime powers can be combined \'{a} la the Chinese Remainder Theorem and some elbow grease to give the Explicit Strong Approximation Theorem below and the argument in \cite{thin} gives the corollary. 
\begin{thm}[Explicit Strong Approximation Theorem] \label{big3} For $(k,6)=1,$ 
\begin{equation*}
    \Gamma(2)'\textit{ $\bmod{(2^m 3^n k)}$ }  \textit{ $\cong$ } \Gamma(2)' \bmod{2^m} \times \Gamma(2)' \bmod{3^n}\times \SL_2(k).
\end{equation*}
\end{thm}
\begin{cor}\label{zariskicoro} The group $\Gamma(2)'$ is Zariski Dense with respect to $\SL_2$.
\end{cor}

\subsection{The Case when $k$ is such that $(k,6)=1$}\label{modprimepowers}

We will note a few lemmas that will help in proving what $\Gamma(2)'\bmod{p^l}$ looks like for $p>3.$ The first lemma is a consequence of Hensel's Lemma and the argument is standard.
\begin{lem}\label{neg}
For $k$ odd, $-I\in \SL_2(k)'$.
\end{lem}

\begin{lem}\label{equiv}
 For $(k,2)=1,$ $\SL_2(k) = \Gamma(2) \bmod{k} $ and $\SL_2(k)' = \Gamma(2)' \bmod{k}$.
\end{lem}
\begin{proof}
First, we note that $\Gamma(2) \leq \SL_2(\mathbb{Z})$. So what remains to be shown is $\SL_2(k) \leq \Gamma(2) \bmod{k}$. Now by assumption $k$ is not divisible by 2. Thus there exists $\alpha$ such that $2\alpha\equiv 1 \bmod{k}$. This easily gives us the next two congruences,  $C\equiv \mathcal{A}^\alpha \bmod{k} \text { and } D\equiv \mathcal{B}^\alpha \bmod{k}.  $ Since $\SL_2(\mathbb{Z})$ is generated by $C, D, -I$ and we just saw that $C, D,-I\in \Gamma(2) \mod k$ we get that $\SL_2(k) = \Gamma(2) \mod k$. Thus $\SL_2(k)'= \Gamma(2)' \mod k$.
\end{proof}

\begin{lem}\label{CD}
For $(k,6)=1,$ $ C,D \in \SL_2(k)'$.
\end{lem}
\begin{proof}
Let $q\equiv -3 \bmod{k}$, $u\equiv q^{-1}\bmod{k}$, $n\equiv 1-u\bmod{k}$, and $m\equiv -2 n^{-1}\bmod{k}$. First, we should justify that $u$ and $m$ actually exist. Given that $k$ is not divisible by 3 we know that $q$ is invertible and thus $u$ exists.  Since $(k,2)=1$, there exists a $t$ such that $4t\equiv 1\bmod{k}$ which gives us $(1-u)3t\equiv 3t+t(-3)u\equiv 3t+t(-3)q^{-1}\equiv 4t\equiv 1 \bmod{k}$ and so $1-u$ is also invertible which gives us that $m$ exists.
Next, we take a specific element in the commutator subgroup and replace $m$ and $u$ as defined above:
\begin{equation*}
  [C^n,D^m][C^u,D^q]\equiv [C^n,D^{-2n^{-1}}][C^{q^{-1}},D^q]\equiv \begin{pmatrix} 9+2nq & -3 q^{-1}\\ 12n^{-1}+3q & -4 (nq)^{-1} \end{pmatrix} \bmod{k}.
\end{equation*}
Modular arithmetic gives us then that $[C^n,D^m][C^u,D^q]\equiv C \bmod{k}$. For the same reason we see that $D=C^T\equiv ([C^n,D^m][C^u,D^q])^T=[C^{-q},D^{-u}][C^{-m},D^{-n}].$ And so $C,D \in \SL_2(k)'.$ 
\end{proof}

\begin{cor}\label{primepowers}
Let $k$ odd and coprime to 3, then $ \SL_2(k) = \SL_2(k)'$.
\end{cor}
\begin{proof}
Given that $C, D, -I$ generate $\SL_2(\mathbb{Z})$ and they are all elements of $\SL_2(k)'$ via Lemmas \ref{neg} and \ref{CD} we have proven this corollary.
\end{proof}

\begin{cor}
 For $k$ odd and coprime to 3,  $\SL_2(k) = \Gamma(2)' \bmod{k}$.
\end{cor}
\begin{proof}
By Lemma \ref{equiv} and Corollary \ref{primepowers}, the statement is proven.
\end{proof}

\subsection{The Case when  $k=2^l$}\label{mod2}

First we will state a lemma that is obvious but will help in the case of $2^l$ and $3^l$.
\begin{lem}\label{helper} Let $p$ be a prime and $n_p \geq 1.$ Then the following two statements are equivalent.
\begin{enumerate}
\item For $l\geq n_p,$  $\Gamma(2)' \bmod{p^l}= \{\gamma \in \SL_2(p^l) : \gamma \in \Gamma(2)' \bmod{p^{n_p}}\}.$
\item For $l\geq n_p,$ $|\Gamma(2)' \bmod{p^{l+1}}|=p^3 \cdot |\Gamma(2)' \bmod{p^{l}}|$
\end{enumerate}

\end{lem}
\begin{proof}
The proof comes naturally after we note two simple facts about the groups involved:

\textbf{\textit{Fact (i)}} For $l\geq n_p,$  $\Gamma(2)' \bmod{p^l}\leq \{\gamma \in \SL_2(p^l) : \gamma \in \Gamma(2)' \bmod{p^{n_p}}\}.$ This is obvious.

\textbf{\textit{Fact (ii)}} For $l\geq n_p,$
        $$\#\{\gamma \in \SL_2(p^{l+1}) : \gamma \in \Gamma(2)' \bmod{p^{n_p}}\}=p^3 \cdot \#\{\gamma \in \SL_2(p^l) : \gamma \in \Gamma(2)' \bmod{p^{n_p}}\}.$$
Indeed, call the left-hand set $X$ and the right-hand set, $Y$. Let $M\in Y.$ Then  $M+p^lM_2 \in X$ if and only if $\det(M+p^lM_2)\equiv 1 \bmod{p^{l+1}}.$ Write  $$M+p^lM_2=\begin{pmatrix}a_1 & b_1\\c_1 & d_1 \end{pmatrix}+p^l \begin{pmatrix}a_2 & b_2\\c_2 & d_2 \end{pmatrix} \bmod{p^{l+1}}.$$ Let $x\in \mathbb{Z}/p\mathbb{Z}$ be such that $\det(M)\equiv 1+p^l x \bmod{p^{l+1}}.$ Then the set of possible values of $M_2$ is a three dimensional vector space over $\mathbb{Z}/p\mathbb{Z}$ which then has size $p^3$ given that $M_2$ must satisfy $x+a_1 d_2+a_2d_1-b_1c_2-b_2c_1\equiv 0 \bmod{p}.$

Proof that $1\Rightarrow 2$
\begin{align*}
    |\Gamma(2)' \bmod{p^{l+1}}|&= \#\{\gamma \in \SL_2(p^{l+1}) : \gamma \in \Gamma(2)' \bmod{p^{n_p}}\} \textit{ by Assuming 1}\\
    &=p^3 \cdot \#\{\gamma \in \SL_2(p^l) : \gamma \in \Gamma(2)' \bmod{p^{n_p}}\} \textit{ by Fact (ii)}\\
    &= p^3 \cdot |\Gamma(2)' \bmod{p^{l}}| \textit{ by Assuming 1 again}
\end{align*} Proof that $2\Rightarrow 1$ This statement is equivalent to showing that the size of the two groups are the same since Fact (i) told us the one is a subgroup of the other. We will prove this via induction. For $l=n_p,$ clearly 1 holds by definition. Now assume that it holds for some $l\geq n_p.$ Then,
\begin{align*}
    |\Gamma(2)' \bmod{p^{l+1}}|&=p^3 \cdot |\Gamma(2)' \bmod{p^{l}}| \textit{ by Assuming 2}\\
    &=p^3 \cdot \#\{\gamma \in \SL_2(p^l) : \gamma \in \Gamma(2)' \bmod{p^{n_p}}\} \textit{ by inductive step}\\
    &=\#\{\gamma \in \SL_2(p^{l+1}) : \gamma \in \Gamma(2)' \bmod{p^{n_p}}\} \textit{ by Fact (ii)}.\end{align*}\end{proof}

Now the goal of this section is to find what $\Gamma(2)' \bmod{2^l}$ is. One can easily calculate by hand that $\Gamma(2)' \bmod{2}= \Gamma(2)' \bmod{4}= \{I\}$, $\Gamma(2)' \bmod{8}= \{I, 5I\}$, and 
\begin{equation*}
  \Gamma(2)' \bmod{2^4}=\left\{I, 9I, \begin{pmatrix}5 & 8\\8 & 13\end{pmatrix},\begin{pmatrix}13 & 8\\8 & 5\end{pmatrix}\right\}.
\end{equation*}

What remains to be proven is that for $l\geq 4,$ case 2 of Lemma \ref{helper} holds. To do this, we define a collection of functions $\left<\rho_{l,i}\mid 1\leq i\leq 3\right>$ such that for any $\rho \in \left<\rho_{l,i}\right>$:
\begin{itemize}
    \item For $M \in \Gamma(2)',$ $\rho(M) \in \Gamma(2)',$
    \item For $M \in \Gamma(2)' \bmod{2^l},$ $\rho(M) \equiv M \bmod{2^l},$
    \item For $M \in \Gamma(2)' \bmod{2^{l+1}},$ $\rho(M) \equiv M + 2^l y_M \bmod{2^{l+1}}$ where we refer to $y_M$ as a shift of $M$.
\end{itemize}

These functions will be specifically chosen so that, for any $M\in \Gamma(2)' \bmod{2^{l+1}},$ all combinations of these functions on $M$ will result in a total of $8=2^3$ unique matrices modulo $2^{l+1}$ that are the same modulo $2^l$ and that all come from a matrix that lies in $\Gamma(2)'.$ In doing this, we will have proven that, for $p=2$ and all $l\geq 4,$ we have the desired result that $|\Gamma(2)' \bmod{p^{l+1}}|=p^3 \cdot |\Gamma(2)' \bmod{p^{l}}|$.

Let $v=2^{l-2}$ and $w=2^{l-4}.$ We will now define this family of functions to be: 
\begin{align*}
\rho_{l,1}(M)&= [\mathcal{A}, \mathcal{B}^v] M\\
\rho_{l,2}(M)&=[\mathcal{A}, \mathcal{B}^{w}][ \mathcal{A}^{-1}, \mathcal{B}^{-3 w}] M\\
\rho_{l,3}(M)&= [\mathcal{A}^{3w},\mathcal{B}] [\mathcal{A}^{-w}, \mathcal{B}^{-1}] M.
\end{align*}

Next, we record how these act on $M\in \Gamma(2)' \bmod{2^{l+1}}$ to demonstrate that they have the desired properties. Write $M\equiv \begin{pmatrix}a_1 & b_1\\c_1 & d_1 \end{pmatrix}+2^l \begin{pmatrix}a_2 & b_2\\c_2 & d_2 \end{pmatrix} \bmod{2^{l+1}}$ where $0\leq a_1, b_1, c_1, d_1 <2^l$ and $0\leq a_2, b_2, c_2, d_2 <2.$ It is also worth mentioning that, since $M$ is in the level-2 principal congruence subgroup, $a_1\equiv d_1\equiv 1 \bmod{2}$ and  $b_1\equiv c_1\equiv 0 \bmod{2}.$

A tedious yet straight forward calculation (that can be seen in the Mathematica folder on the GitHub repository) will give us the following equivalences (where $\delta_l$ is 0 if $l=4$ and 1 if $l>4$):
  \begin{equation*}
      \rho_{l,1}(M)\equiv M+2^l \begin{pmatrix}1 & 0\\0 & 1 \end{pmatrix} \bmod{2^{l+1}}
  \end{equation*}

  \begin{equation*}
    \rho_{l,2}(M)\equiv M+2^l \begin{pmatrix} \delta_l & 1\\0 & \delta_l \end{pmatrix}  \bmod{2^{l+1}}
  \end{equation*}
 \begin{equation*}
    \rho_{l,3}(M)\equiv M+2^l \begin{pmatrix} \delta_l & 0\\1 & \delta_l \end{pmatrix} \bmod{2^{l+1}}.
  \end{equation*} Note that the shift matrices all lie in the lie algebra of $\SL_2(\mathbb{Z}/2\mathbb{Z})$, that is the trace of the shifts are congruent to $0 \bmod{2}.$ Thus the shifts, for $l\geq 4,$ are a 3-dimensional vector space over $\mathbb{Z}/2\mathbb{Z}$ and have cardinality $2^3=8.$ Combining this with Lemma \ref{helper} we get the following.
\begin{cor}
For $l\geq 4,$ $\Gamma(2)' \bmod{2^l}= \{\gamma \in \SL_2(2^l) : \gamma \in \Gamma(2)' \bmod{2^4}\}$.
\end{cor}

\subsection{The Case when $k=3^l$}\label{mod3}

This case will be done almost identically to the case where $k=2^l.$ However, we will utilize the fact that  $\gcd(k, 2)=1,$ which, by Lemma \ref{equiv}, allows us to answer the question for $\SL_2(k)'$.  We wish to prove that for $l\geq 1,$ case 2 of Lemma \ref{helper} holds. To do this, we form a family of functions $\left<\sigma_{l,i}\mid 1\leq i\leq 6\right>$ such that for $\sigma \in \left<\sigma_{l,i}\right>$:
\begin{itemize}
    \item For $M \in \SL_2(\mathbb{Z})',$ $\sigma(M) \in \SL_2(\mathbb{Z})',$
    \item For $M \in \SL_2(3^l)',$ $\sigma(M) \equiv M \bmod{3^l},$
    \item For $M \in \SL_2(3^{l+1})',$ $\sigma(M) \equiv M + 3^l y_M \bmod{3^{l+1}}$ where $y_M$ is a shift.
\end{itemize}

Again, the idea is to define these functions such that for all $M\in \SL_2(k)',$ the family of functions applied to $M$ will result in $27=3^3$ unique matrices modulo $3^{l+1}$ that are equivalent modulo $3^l.$ Then we will have proven that, for all $l\geq 1,$ it is true that $|\SL_2(3^{l+1})'|=3^3\cdot |\SL_2(3^l)|.$

Now let $t_1:=C^{3^l}$ and $t_2:=D^{3^l}.$ Then the family of functions will be the following.
\begin{align*}
 \sigma_{l,1}(M)= t_1 M t_1^{-1} &&
\sigma_{l,2}(M)=C^{-1} t_2 C M t_2^{-1} &&
\sigma_{l,3}(M)=C t_2 C^{-1} M t_2^{-1}\\
\sigma_{l,4}(M)=t_2 M C^{-1} t_2^{-1} C &&
\sigma_{l,5}(M)=t_1 M D^{-1} t_1^{-1} D &&
\sigma_{l,6}(M)=D^{-1} t_1 D M t_1^{-1}
\end{align*}

 Again, we record how these functions act on $M\in \SL_2(3^{l+1})$ and mention that the Mathematica file for this section verifies these statements. Write 
 \begin{equation*}
     M\equiv \begin{pmatrix}a_1 & b_1\\c_1 & d_1 \end{pmatrix}+3^l\begin{pmatrix}a_2 & b_2\\c_2 & d_2 \end{pmatrix} \bmod{3^{l+1}}
 \end{equation*} where $0\leq a_1, b_1, c_1, d_1 <3^l$ and $0\leq a_2, b_2, c_2, d_2 <3.$ Then we see that these operations are:
\begin{align*}
 \sigma_{l,1}(M)&\equiv M+3^l \begin{pmatrix}c_1& d_1-a_1\\0& -c_1 \end{pmatrix} \bmod{3^{l+1}} \\
 \sigma_{l,2}(M)&\equiv M+3^l \begin{pmatrix}-a_1-b_1-c_1& -b_1-d_1\\a_1+c_1-d_1& b_1+d_1 \end{pmatrix} \bmod{3^{l+1}} \\
\sigma_{l,3}(M)&\equiv M+3^l \begin{pmatrix}a_1-b_1-c_1& b_1-d_1\\a_1-c_1-d_1& b_1-d_1 \end{pmatrix} \bmod{3^{l+1}}\\
\sigma_{l,4}(M)&\equiv M+3^l \begin{pmatrix}a_1-b_1& a_1-b_1\\a_1+c_1-d_1& b_1+c_1-d_1 \end{pmatrix} \bmod{3^{l+1}}\\
\sigma_{l,5}(M)&\equiv M+3^l \begin{pmatrix}-a_1+b_1+c_1& -a_1+b_1+d_1\\-c_1+d_1& -c_1+d_1 \end{pmatrix}\bmod{3^{l+1}} \\
\sigma_{l,6}(M)&\equiv M+3^l \begin{pmatrix}a_1+c_1& -a_1+b_1+d_1\\-a_1-c_1& -b_1-c_1-d_1 \end{pmatrix}\bmod{3^{l+1}}.
\end{align*}

\begin{lem}\label{linearCombos} The set of all linear combinations of these shifts (for which there are 729) evaluated at any $M\in M_2(\mathbb{Z}/3\mathbb{Z})$ that is not the zero matrix, gives exactly 27 unique shifts mod 3.
\end{lem} The proof of this lemma is via direct computation and can be found in the Mathematica file.

\begin{cor}
For $l\geq 1,$ $\SL_2(3^l)'= \{\gamma \in \SL_2(3^l) : \gamma \in \SL_2(3)'\}$. 
\end{cor}

\begin{proof}
By Lemma \ref{linearCombos} above, we know that for each $M \in \SL_2(3^{l+1})'$ there are a combination of actions on $M$ that will result in 27 different shifts i.e. 27 different matrices modulo $3^{l+1}$ which are all the same modulo $3^l$ and so by Lemma \ref{helper} we are done.
\end{proof}

\begin{cor}
For $l\geq 1,$ $\Gamma(2)' \bmod{3^l}= \{\gamma \in \SL_2(3^l) : \gamma \in \Gamma(2)' \bmod{3}\}$. 
\end{cor}
\begin{proof}
By Lemma \ref{equiv} and the previous corollary, the statement is proven.
\end{proof}

\subsection{The Case of general $k$}

First, we will state and some general group theoretic properties (that together will be similar to the classic Goursat's Lemma) which will help in piecing the previous three portions of this section together. For the remainder of this section, given a subgroup of $G\times K$, $\phi_G$ is a projection onto $G$ while $\phi_K$ is a projection onto $K$.
\begin{lem}\label{norm}
Let $N\lhd G\times K$ and let the projection of $N$ onto $K$ be surjective, then $(e, [K,K])\leq N,$ where $e$ is the identity of $G.$
\end{lem}

\begin{lem}\label{gora} 
Let $N\lhd G\times K$, $\phi_G(N)=G$ and $ker(\phi_G)\cong K$. Then $N=G\times K.$
\end{lem}

\begin{lem}\label{almost} Let $(k,6)=1$, 
\begin{equation*}
  \Gamma(2)'\bmod{2^m 3^n k}\cong \Gamma(2)' \bmod{2^m 3^n} \times \SL_2(k).
\end{equation*} 
\end{lem}

\begin{proof}
Let
\begin{align*}
  N&=\{(\gamma\bmod{2^m3^n},\gamma \bmod{k}):\gamma \in \Gamma(2)' \bmod{ 2^m 3^n k}\}\\
  G&=\Gamma(2)' \bmod{2^m 3^n}\\
  K&=\Gamma(2)' \bmod{k}=\SL_2(k).
\end{align*}

Recall that $\Gamma(2)'$ is a normal subgroup of $\SL_2(\mathbb{Z})$. Now take any $(g,h) \in G\times K.$ Then by the Chinese Remainder Theorem, there exists $\alpha \in \SL_2(\mathbb{Z})$ such that $g \equiv \alpha \bmod 2^n 3^m$ and $h \equiv \alpha \bmod k.$ And so $N\lhd G\times K$ given that for any $(\gamma\bmod{2^m3^n},\gamma \bmod{k})\in N$, we see that 
\begin{equation*}
  (g,h)(\gamma\bmod{2^m3^n},\gamma \bmod{k})(g^{-1}, h^{-1})\equiv (\alpha \gamma \alpha^{-1} \bmod{2^m3^n}, \alpha \gamma \alpha^{-1} \bmod{k})\in N.
\end{equation*} 

Now let $\phi_X(N)$ for $X\in \{G, K\}$ be the projection of $N$ onto $X$. Take $g\in G,$ then there exists $\gamma \in \Gamma(2)'$ such that $\gamma \equiv g\bmod 2^m 3^n$ and $(\gamma \bmod 2^m 3^n, \gamma \bmod k)\in N$. Thus $\phi_G(N)=G.$ The same reasoning gives us that $\phi_K(N) =K.$ Finally, combining Lemmas \ref{norm} and \ref{gora}, we have proven the claim.
\end{proof}

The goal is now to prove Theorem \ref{big3} which, given Lemma \ref{almost}, reduces to just proving that $ \Gamma(2)'\bmod{2^m 3^n }\cong \Gamma(2)' \bmod{2^m} \times \Gamma(2)' \bmod{3^n}.$  By definition, we have $$\Gamma(2)'\bmod{2^n 3^m }\leq \Gamma(2)' \bmod{2^n} \times \Gamma(2)' \bmod{3^m}$$ and when $n=0$ it's true for all $m$. 

The idea for the remainder of the proof is to first modify the family of functions $\left<\rho_{l,i}\right>$ constructed in section \ref{mod2} in a way that doesn't change the original set of properties but is more careful in the way that it affects elements modulo 3, then we do a similar change in $\left< \sigma_{l,i}\right>$ constructed in Section \ref{mod3} in a way that, again, doesn't change the original set of properties but is careful in the way that it affects elements modulo $2^m$.

\begin{lem} Let $m\geq 0$ then
\begin{equation}\label{nm1}
  \Gamma(2)'\bmod{2^m 3^1}\cong \Gamma(2)' \bmod{2^m} \times \Gamma(2)' \bmod{3}.
\end{equation} 
\end{lem}
\begin{proof}
First, via a direct computation, we can see that, for the finite number of cases where $0\leq m\leq 4,$ (\ref{nm1}) holds. Suppose that the equation holds for some $m\geq 4.$ Let us then change the actions $\left<\rho_{m,i}\right>$ from Section \ref{mod2} to have the first three conditions as before but also be such that $\rho_{m,i}(M)\equiv M \bmod{3}.$

A simple calculation will verify that the new family of functions,
\begin{equation*}
    \rho_{m,i}^*(M)=\begin{cases}[\mathcal{A}^{2^m}, \mathcal{B}]\rho_{m,i}(M) & \textit{ if $i=3$ and $m$ is even}\\
    [\mathcal{B}^{2^m}, \mathcal{A}]\rho_{m,i}(M) & \textit{ else,}
    \end{cases}
\end{equation*} satisfies these conditions. So for $M\in \Gamma(2)' \bmod{2^{m+1} 3^1}$ the family of $\left<\rho_{m,i}^*\right>$ gives exactly 8 unique matrices in  $\Gamma(2)' \bmod{2^{m+1}3^1}$ that are the same modulo $2^m 3^1.$ Thus,
\begin{equation*}
    |\Gamma(2)' \bmod{2^{m+1} 3^1}|= 8 \cdot |\Gamma(2)' \bmod{2^m 3^1}|=|\Gamma(2)' \bmod{2^{m+1}}\times \Gamma(2)' \bmod{3}|.
\end{equation*} And this  completes the proof.
\end{proof}

\begin{lem}\label{mn} Let $m,n\geq 0,$ then
\begin{equation}\label{nm}
  \Gamma(2)'\bmod{2^m 3^n }\cong \Gamma(2)' \bmod{2^m} \times \Gamma(2)' \bmod{3^n}.
\end{equation} 
\end{lem}
\begin{proof}
Fix $m$. Then we know that the statements hold for this $m$ and $n=1$ by the previous lemma. Now suppose it holds true for $m$ and some $n\geq 1.$  Then we will prove that it must hold for this $m$ and $n+1$ in the following way. Recall that the any action of $\left<\sigma_{n,i}\right>$ is just a combination of powers of $C$ and $D.$ However when we say $C$ we are hiding the fact that we actually mean a power of $\mathcal{A}$ that is congruent to $C$ modulo $3^{n+1}$ and that $D$ is that same power of $\mathcal{B}$ that is congruent modulo $3^{n+1}.$ In Lemma \ref{CD} this power was described as $\alpha$ such that $2\alpha \equiv 1 \bmod{3^{n+1}}$. However, there are infinitely many such $\alpha$ and we can choose this $\alpha$ to have the added property that $\alpha \equiv 0 \bmod{2^m}$ given that we fixed $m$. This property means that for any $M\in \Gamma(2)' \bmod{2^m 3^{n+1}}$, $\sigma_{n,i}(M)\equiv M \bmod{2^m}$ and so  this family of functions is such that for any such $M$ there are 27 unique matrices modulo $2^m 3^{n+1}$ and that are all the same modulo $2^m 3^n.$ Which completes the proof since 
\begin{equation*}
    |\Gamma(2)' \bmod{2^m 3^{n+1}}|= 27 \cdot |\Gamma(2)' \bmod{2^m 3^{n}}|=|\Gamma(2)' \bmod{2^m}\times \Gamma(2)' \bmod{ 3^{n+1}}|.
\end{equation*}
\end{proof}

Now we may finally conclude the Strong Approximation Theorem \ref{big3} by  plugging Lemma \ref{mn} into Lemma \ref{almost}.


\section{Calculating the Local Traces}\label{ch:Trace}

When one looks at the representatives of the groups generated in the previous section, certain theorems about the traces in the local case are hinted at and we shall spend this section proving them. This section relies on work in the GitHub repository. \footnote{\href{https://github.com/BrookeOgrodnik/CommutatorSubgroups/blob/master/Mathematica/tracesmodk.nb}{https://github.com/BrookeOgrodnik/CommutatorSubgroups/blob/master/Mathematica/\newline tracesmodk.nb}}
\begin{lem}\label{sltr}
When $(k,6)=1$, $\Tr(\Gamma(2)') \bmod{k}=\{0,1,\ldots, k-1\}$.
\end{lem}
\begin{proof}
Given $(k,6)=1,$ we know from Theorem \ref{big} that $\Gamma(2)' \bmod{k}= \SL_2(k)$ and since, for $0\leq a< k$, $\begin{pmatrix}0 & -1 \\ 1 & a\end{pmatrix}\in \SL_2(k),$ we conclude all traces appear modulo $k.$
\end{proof}

\begin{lem}\label{mod3tr} \begin{equation*}
   \Tr(\Gamma(2)') \bmod 3^m =\begin{cases}\{0,1,2\} &\text{if } m=1\\
    \{0,2,3,6,7\} &\text{if } m=2\\
    \{ 0\leq t<3^m \mid t \in \Tr(\Gamma(2)') \bmod 9\} &\text{if } m>2.
    \end{cases}
\end{equation*}
\end{lem}

\begin{remark}
For any set $X,$ if $|X \bmod{p^m}|=n$ then $|X \bmod{p^{m+1}}|\leq pn.$ So if, for each $m\geq 2$, we can create a subset $X_m^{(3)}$ of elements in $\Gamma(2)'$ that has the property that $X_{m+1}^{(3)}\cong X_m^{(3)} \bmod{3^m}$ and $|\Tr(X_{m+1}^{(3)}) \bmod{3^{m+1}}|= 3 \cdot | \Tr(X_{m}^{(3)}) \bmod{3^{m}}|,$ then we will have proven this lemma since the maximum number of elements appear. This note also applies for later when we prove a similar lemma for the powers of 2. 
\end{remark}

\begin{proof}
In the associated Mathematica file, this theorem is first verified via brute force for the cases where $m\leq 2.$ It is here we see that $\Tr(\Gamma(2)') \bmod 9= \{0, 2, 3, 6, 7\}.$

Recall the action of $\sigma_{m,3}(M).$ As before, we write \begin{equation*}
     M\equiv \begin{pmatrix}a_1 & b_1\\c_1 & d_1 \end{pmatrix}+3^m\begin{pmatrix}a_2 & b_2\\c_2 & d_2 \end{pmatrix} \bmod{3^{m+1}}
 \end{equation*} where $0\leq a_1, b_1, c_1, d_1 <3^m$ and $0\leq a_2, b_2, c_2, d_2 <3$ and from here we can write out how $\sigma_{m,3}$ acts on $M$ modulo $3^{m+1}$ and modulo $3^{m+2}.$
 \begin{align*}
     \textit{For $m\geq 1, $  }\sigma_{m,3}(M)&\equiv M+3^m \begin{pmatrix}a_1-b_1-c_1& b_1-d_1\\a_1-c_1-d_1& b_1-d_1 \end{pmatrix} \bmod{3^{m+1}}.\\
      \textit{For $m\geq 2, $  }\sigma_{m,3}(M)&\equiv M+3^m \begin{pmatrix}a_1-b_1-c_1& b_1-d_1\\a_1-c_1-d_1& b_1-d_1 \end{pmatrix} \bmod{3^{m+2}} .
 \end{align*} Note here that in order for us to be able to look modulo $3^{m+2}$ in general, we need $m\geq 2$ and so we will need to handle the cases when $m=0,1$ separately (this fact will pop up again when proving the similar lemma for powers of 2). We will also comment that it is not an option to look modulo $3$ for the cases where $\Tr(M) \bmod{9}$ is in $\{2,7\}$ because for such $M$'s and all $\sigma\in \left<\sigma_{m,i}\right>,$ $\Tr(\sigma(M))\equiv \Tr(M) \bmod{3^{m+1}}$ which is why we must look modulo $3^{m+2}.$
 
In both cases one can easily calculate the trace.
 \begin{align*}
     \textit{For $m\geq 1, $  }\Tr(\sigma_{m,3}(M))&\equiv \Tr(M)+3^m (a_1-c_1-d_1) \bmod{3^{m+1}}.\\
     \textit{For $m\geq 2, $  }\Tr(\sigma_{m,3}(M))&\equiv \Tr(M)+3^m (a_1-c_1-d_1) \bmod{3^{m+2}}.
 \end{align*}
 
 So what we can see is that the action of $\sigma_{m,3}$ on $M$ only depends on what $M$ is modulo 3 in the first case and modulo 9 in the second. 

Define $M_1:=[\mathcal{A},\mathcal{B}],$ $M_2:=[\mathcal{A}^5,\mathcal{B}^{-6}],$ and $M_3:=M_2^{-1}M_1^2.$ Now we can go on to define the subsets of $\Gamma(2)'$.

For $m\geq 1,$ let $\chi_{1}(M_1):=\{M_1\}$ and
\begin{equation*}
    \chi_{m+1}(M_1):=\bigsqcup_{M\in \chi_{m}(M_1)}\{M, \sigma_{m,3}(M), \sigma_{m,3}^2(M)\}.
\end{equation*}

For $i=2,3$, let $ \chi_{2}(M_i):=\{M_i\}$ and
\begin{equation*}
     \chi_{3}(M_i):=\begin{cases}
     \{M_2, \sigma_{1,3}(M_2), \sigma_{1,5}(M_2)\} &\textit{ if $i=2$}\\
      \{M_3, \sigma_{1,2}(M_3), \sigma_{1,5}(M_3)\} &\textit{ if $i=3$}.
     \end{cases} 
\end{equation*} Then for $m\geq 2$ let
\begin{equation*}
      \chi_{m+2}(M_i):=\bigsqcup_{M\in \chi_{m}(M_i)}\{M, \sigma_{m,3}(M), \sigma_{m,3}^2(M)\}.
\end{equation*} 
Finally, for $m\geq 2,$
\begin{equation*}
     X_m^{(3)}:=\bigsqcup_{1\leq i\leq 3} \chi_{m}^{(3)}(M_i) \text{ and }tr_m^{(3)}:=\Tr(X_m^{(3)}) \bmod{3^m}.
\end{equation*}
We will show that these subsets of $\Gamma(2)'$ have the desired properties and we will show that $tr_m^{(3)}=\Tr(\Gamma(2)') \bmod{3^m}$.

For $m\geq 1,$ $\Tr(\chi_{1}(M_1))=\{0\} \bmod{3}$ and 
\begin{equation*}
    \Tr(\chi_{m+1}(M_1))=\bigsqcup_{M\in \chi_{m}(M_1)}\{\Tr(M), \Tr(M)+3^m, \Tr(M)+2\cdot 3^m\} \bmod{3^{m+1}}.
\end{equation*} 

In the next case we have $ \Tr(\chi_{2}(M_2))=\{2\} \bmod{9},$ $\Tr(\chi_{3}(M_2))=\{ 2,11, 20\} \bmod{27},$ and, for $m\geq 2,$
\begin{equation*}
    \Tr(\chi_{m+2}(M_2))= \bigsqcup_{M\in \chi_{m}(M_2)}\{\Tr(M), \Tr(M)+3^{m+1}, \Tr(M)+2\cdot 3^{m+1}\} \bmod{3^{m+2}}.
\end{equation*} And here, $ \Tr(\chi_{2}(M_3))=\{7\} \bmod{9},$ $\Tr(\chi_{3}(M_3))=\{7, 16, 25\} \bmod{27},$ and, for $m\geq 2,$
\begin{equation*}
    \Tr(\chi_{m+2}(M_3))= \bigsqcup_{M\in \chi_{m}(M_3)}\{\Tr(M), \Tr(M)+3^{m+1}, \Tr(M)+2\cdot 3^{m+1}\} \bmod{3^{m+2}}.
\end{equation*}

Combining these statements gives $tr_2^{(3)}= \{0, 2, 3, 6, 7\}$, which matches the trace set for $\Gamma(2)' \bmod{9},$ and, for $m\geq 3$, $$ |tr_{m+1}^{(3)}|= 3 \cdot | tr_{m}^{(3)}|$$ So this construction not only gives us the lemma for $m\geq 2,$ but is also a way that, given any trace in $ tr_{m}^{(3)},$ we can easily construct a (probably very large) element with that trace modulo $3^{m}$ via this construction explained above.
\end{proof}

The proof for the case of powers of two case will be similar to the previous lemma.

\begin{lem}\label{mod2tr} \begin{equation*}
   \Tr(\Gamma(2)') \bmod 2^m =\begin{cases}\{0\} &\text{if } m=1\\
    \{2\} &\text{if } m=2,3,4\\
    \{2,18\} &\text{if } m=5,6\\
    \{2,18,66\} &\text{if } m=7\\
    \{2,18,66, 146\} &\text{if } m=8\\
    \{ 0\leq t<2^m \mid t \in \Tr(\Gamma(2)') \bmod 2^8\} &\text{if } m>8.
    \end{cases}
\end{equation*}
\end{lem}

\begin{proof}
Keep in mind that we can (and did in the Mathematica file) verify the statement for $m\leq 8.$ Now, as in the case of powers of 3, we choose functions that have the correct traces modulo $2^8$ and then show that the sizes double after this point as we raise $m.$ Define $M_1:=[\mathcal{A},\mathcal{B}^4 ][\mathcal{A}^{3},\mathcal{B}^{1}] [\mathcal{A}^{-1}, \mathcal{B}^{-1}]$, $M_2:=[\mathcal{A},\mathcal{B}]$, $M_3:=[\mathcal{A},\mathcal{B}^2]$, and $M_4:=[\mathcal{A},\mathcal{B}^5]^{-1}.$ Then let $X_8^{(2)}=\{M_1, M_2, M_3, M_4\}$ and $t_m^{(2)}:=\Tr(X_m^{(2)})\bmod{2^m}$ for $m\geq 8.$ Again, once we define the rest of the $X_m$'s and show that $|t_{m+1}^{(2)}|=2\cdot |t_m^{(2)}|$ we will have completed this proof. 

If we were to write out how $\rho_{m-5,i}$ acts on $M$ modulo $2^{m}$ we would see that the actions have a dependence on $m$ and this would only be removed when $m>13.$ For this reason we must handle the cases between $9\leq m\leq 13$ independently from the cases where $m>13.$ Because there are so many cases, the Mathematica file contains the work in detail but here we just give the punchline via  Table \ref{actions}.

\begin{table}[h]
\begin{center}
{\renewcommand{\arraystretch}{1.25}
 \begin{tabular}{||c|| c | c | c | c |} 
 \hline 
 $m\backslash M$ & $M_1$ & $M_2$ & $M_3$ & $M_4$ \\ [0.5ex] 
 \hline\hline
 9 & $\rho_{4,2}(M)$ & $\rho_{4,1}^{(2)}(M)$& $\rho_{4,3}(M)$ & $\rho_{4,1}^{(2)}(M)$ \\
 \hline
 10 & $\rho_{5,2}(M)$ & $\rho_{5,2}(M)$& $\rho_{5,2}(M)$ & $\rho_{5,2}(M)$ \\
 \hline
 $>10$ & $\rho_{m-5,2}(M)$ & $\rho_{m-5,3}(M)$ & $\rho_{m-5,2}(M)$ & $\rho_{m-5,3}(M)$ \\
 [1ex] 
 \hline
\end{tabular}}
\caption{The actions on $M$ such that $\Tr(\rho(M))-\Tr(M) = 2^{m-1}$ modulo $2^m$}
\label{actions}
\end{center}
\end{table}
Now we can use this table to define $X_m^{(2)}$ for $m\geq 9.$ Let $F(m,M)$ be the entry in row $m$ of Table \ref{actions} and column $M_i$ (where $M\equiv M_i \bmod 2^5$). Then, for $m\geq 9,$
\begin{equation*}
    X_m^{(2)}=\bigsqcup_{M\in X_{m-1}^{(2)}} \{M, F(m,M)\}.
\end{equation*}
This table was constructed so that $\Tr(F(m,M))\equiv \Tr(M)+2^{m-1} \bmod{2^m}$ which then gives us that $|t_{m+1}^{(2)}|=2\cdot |t_m^{(2)}|$ and we have a constructive algorithm for each local trace.
\end{proof}

Via Lemmas \ref{sltr}, \ref{mod3tr}, and \ref{mod2tr} we get the following.

\begin{cor}
  The set of admissible values for $\Gamma(2)'$ are all $t$ such that 
  \begin{equation*}
      t \bmod{256} \in \{2,18,66,146\} \text{ and  } t \bmod{9} \in \{0,2,3,6,7\},
  \end{equation*} and so the list of these admissible values modulo $9\cdot 256$ are
  \begin{align*}
  \{&2,18,66,146,258,322, 402,578,786,834,1026,1042,\\&1170,1298,1554,1602,1794,1938,2050,2194\}.
\end{align*}
\end{cor}




\section{Existence of Admissible Traces}\label{ch:existence}

Section \ref{ch:Admi} told us which $t$'s are admissible. The next logical question is, do all admissible $t$'s appear as traces in $\Gamma(2)'$. As before, this section relies on work in the GitHub repository. \footnote{\href{https://github.com/BrookeOgrodnik/CommutatorSubgroups/blob/master/Mathematica/class\_number\_search.nb}{https://github.com/BrookeOgrodnik/CommutatorSubgroups/blob/master/Mathematica/ \newline class\_number\_search.nb}}

Certain traces obviously appear.
\begin{lem}\label{lem2}
If $t=16a^2+2$ for some integer $a$, then $\Tr\left(\left[\begin{pmatrix}1 & 2\\ 0 &1 \end{pmatrix}, \begin{pmatrix} 1 & 0 \\ 2a & 1 \end{pmatrix}\right]\right)=t$ and so there exists a matrix $g\in \Gamma(2)'$ such that $\Tr(g)=t.$
\end{lem}

We also see that if a $t$ appears as a trace, then there are many other matrices in $\Gamma(2)'$ that will also have this trace due to the cyclic property of traces. Given that for $g,h\in \SL_2(\mathbb{Z}),$ 
\begin{align*}
  \Tr[g,h]&=\Tr[h,g]=\Tr[{}^t h^{-1},{}^t g^{-1}]=\Tr[h,g^{-1}]\\
  &=\Tr[g^{-1},h^{-1}]=\Tr[h^{-1},g]=\Tr[\pm g, \pm h].
\end{align*}

However, in the case of our question of getting all admissible traces, the answer is no.
\begin{lem}
The only admissible $t$'s such that $|t|<10^6$ and $t$ is not the trace of an element in $\Gamma(2)'$ are when $t\in \{-1006, -1726, -2558, -4718\}.$
\end{lem}

\subsection{Algorithm for Existence}\label{algoExistence}

We will first briefly recall some definitions and properties of Binary Quadratic Forms so that we might use the correspondence between them and hyperbolic matrices to find representatives of matrices with trace $t$ from each $\SL_2(\mathbb{Z})$ conjugacy class. These definitions come mainly from \textit{Applications of Thin Orbits} \cite{AOTO} and Section 8 of \textit{Zetafunktionen und quadratische {K}\"{o}rper} \cite{Zag}.

A binary quadratic form is $f(x,y)=Ax^2+Bxy+Cy^2$ with $A,B,C\in \mathbb{Z}$ which we can also write in shorthand notation as $[A,B,C].$ A binary quadratic form's \textbf{discriminant} is $D=B^2-4A C$.

\begin{defn}
Two binary quadratic forms $f$ and $g$ are \textbf{(narrowly) equivalent} if there exists $a,b,c,d\in \mathbb{Z}$ such that $ad-bc=1$ and $g(x,y)=f(a x+b y,c x+d y).$ 
\end{defn}

A root of a binary quadratic form $f(x,1)=Ax^2+Bx+C$ with $A\neq 0$ and discriminant $D$ is $\theta_f=\frac{-B+\sqrt{D}}{2A}$ when $A\neq 0$. We look at the continued fraction expansion of $\theta_f$ such that $-1<\overline{\theta_f}<0$ and $1<\theta_f$, the continued fraction expansion is
\begin{equation*}
  [\overline{a_1,\ldots, a_k}]=a_1+\frac{1}{a_2+\frac{1}{\cdots+\frac{1}{a_k+\frac{1}{\theta_f}}}}=\theta_f.
\end{equation*}

 Now let $f$ and $g$ be (narrowly) equivalent binary quadratic forms and $\theta_f$ and $\theta_{g}$ be their roots. Then we know there exists $
\gamma=\begin{pmatrix}a & b\\ c &d \end{pmatrix}\in \SL_2(\mathbb{Z})$ such that $f(a x+b y, c x+d y)=g(x,y).$ Note that $\gamma^{-1} \theta_f=\theta_{g}.$
 It is for this reason that we many times judge if two binary quadratic forms are (narrowly) equivalent via looking at their continued fraction expansion and, furthermore, we just need to know if the continued fractional expansion of $\theta_f$ is that of $\theta_{g}$ after some even number of shifts i.e. are equivalent with respect to $\SL_2(\mathbb{Z}).$

\begin{defn}
  We call a binary quadratic form $[A,B,C]$ \textbf{primitive} if $\gcd(A,B,C)=1.$
  \end{defn}
  \begin{defn}
  We define \textbf{the narrow class number of $D$}, $h_n(D)$, to be the number of (narrowly) inequivalent primitive binary quadratic forms of discriminant $D$.
  \end{defn}

   One can use Gauss's reduction algorithm to construct representatives of the equivalency classes of binary quadratic forms. However, this algorithm may find multiple representatives of the same class. So this algorithm gives us a finite list of binary quatric forms but it is only an upper bound for the (narrow) class number. Thus, there is one last step to ensure that the binary quadratic forms are in different classes. This is done via computing the continued fraction expansion of $\theta_f$ where $f$ ranges over the different forms found in Gauss's Algorithm. If one wanted to find the (narrow) class number, they would count the remaining forms that are also primitive. We then use this class number, $h_n(D)$, to answer questions about the traces in $L$ and $\Gamma(2)'$.

\begin{lem}\label{onetoone}
Let $t\neq 2$ be congruent to $2\bmod{16}$ and $D=(t^2-4)/64$. Then there is a bijection between binary quadratic forms of $[A,B,C]$ with discriminant $D$ and hyperbolic matrices with fixed trace $t$ in $L$.
\end{lem}
\begin{proof}
Let $t\equiv 2 \bmod{16}$ and $D>0$ such that $t^2-64D=4.$ Now take, $[A,B,C]$, a binary quadratic form (not necessarily primitive) of discriminant $D$.  Then we can define the following map:
\begin{equation*}
 [A,B,C] \mapsto \begin{pmatrix}
  \frac{t-8B}{2} & -8 C\\ 8A & \frac{t+8B}{2}
  \end{pmatrix}\in L.
\end{equation*}
If $\begin{pmatrix}4a+1 & 8b\\ 8c& 4d+1 \end{pmatrix}\in L,$ then $
[c,(d-a)/2,-b]\mapsto \begin{pmatrix}
\frac{t-4(d-a)}{2} & 8b \\ 8c & \frac{t+4(d-a)}{2}
\end{pmatrix}. $ However note that,
\begin{align*}
  D&=\frac{(d-a)^2}{4}+4bc=\frac{(4(a+d)+2)^2-4}{64}.
\end{align*}
Thus $t=4(a+d)+2$ which means that $\frac{t-4(d-a)}{2}=4a+1$ and $\frac{t+4(d-a)}{2}=4d+1$ and so the defined mapping is onto.
\end{proof}

This lemma is key in helping us to prove a lower bound on the class number, $h(t),$ [see Corollary \ref{corexist}, \cite{dissertation} page 39].  It is also what allows us to calculate $h(t)$ and $h'(t)$ explicitly as we are about to see.

So how exactly do we do this? Exploiting the correspondence between binary quadratic forms and hyperbolic matrices, we can find these representatives systematically. Once we have a set of the unique hyperbolic matrices, we have the $L$ class number of $t$, $h(t)$. Now we want to check which of these matrices, if any, are in $\Gamma(2)'.$ The condition is, of course, that an element is in $\Gamma(2)'$ if and only if the sum of the exponents of $\mathcal{A}$'s is 0 as well as those of $\mathcal{B}$'s, i.e. checking that the matrix is in the trivial homology.

Writing the word in terms of its generators is done simply by using the fundamental domain of $H$ and looking at the matrix's action on $i$. We recall that the closure of the fundamental domain of $H$ is the set, $\{x+iy\in \mathbb{H} \mid -1\leq x\leq 1, (x-1/2)^2+y^2\geq 1/4 \text{ and } (x+1/2)^2+y^2\geq 1/4\}.$

\begin{figure}[h]
    \centering
    \includegraphics[width=0.4\textwidth]{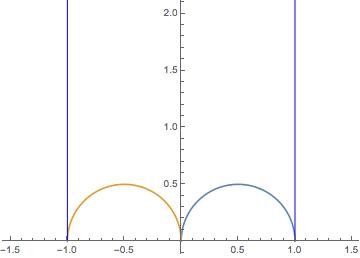}
    \caption{The fundamental domain of $H$}
    \label{gamma2}
\end{figure}

Written out in more detail, let $M$ be the matrix in question. Define $\alpha=gi.$ When $\alpha$ is not in the fundamental domain of $H$ we do the following. If $\alpha$ is such that
\begin{equation*}
  (Re(\alpha)-1/2)^2+Im(\alpha)^2<1/4 \text{ or } (Re(\alpha)+1/2)^2+Im(\alpha)^2<1/4,
\end{equation*} let $n$ be the closest integer to $\frac{-Re(\alpha)}{2(Re(\alpha)^2+Im(\alpha)^2)}$. Then $K=\mathcal{B}^n.$ Else, $\alpha$ has real part either greater than 1 or less than -1 and not in one of the lower two semi-circles. So let $K=\mathcal{A}^n$ where $n$ is such that $-1<2n+Re(\alpha)<1.$  Reassign $\alpha\rightsquigarrow K \alpha$ and repeat until $\alpha$ is back in the fundamental domain of $H$.

\begin{figure}[h!]
\center
\includegraphics[width=60mm]{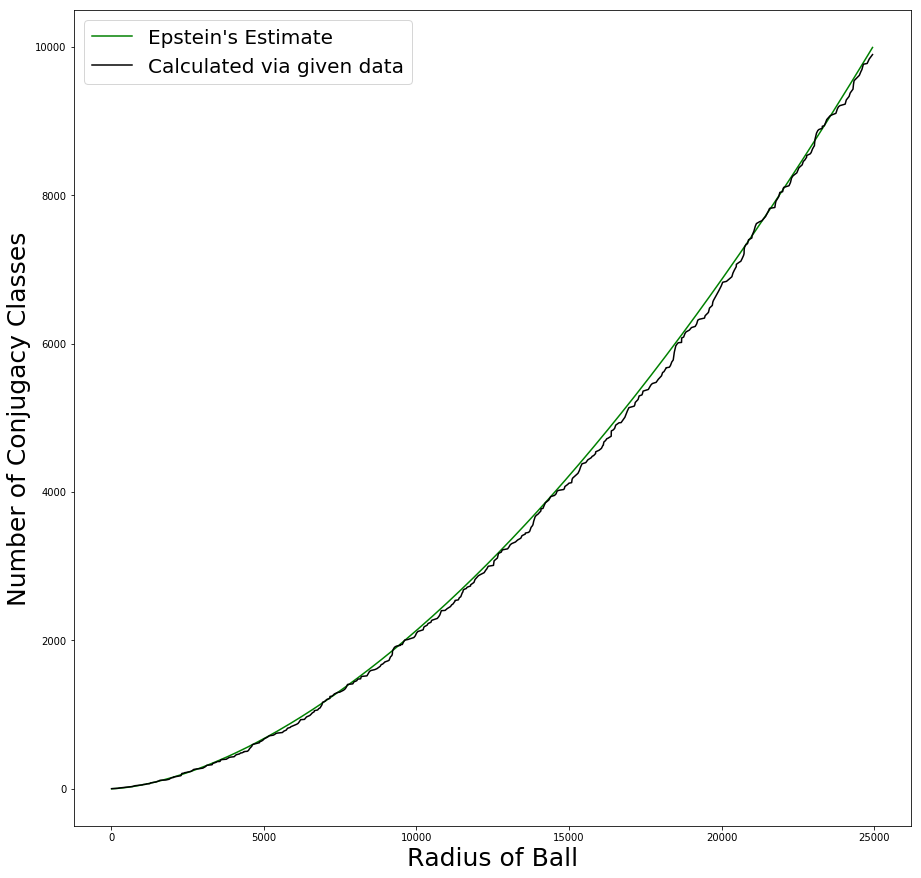}\includegraphics[width=60mm]{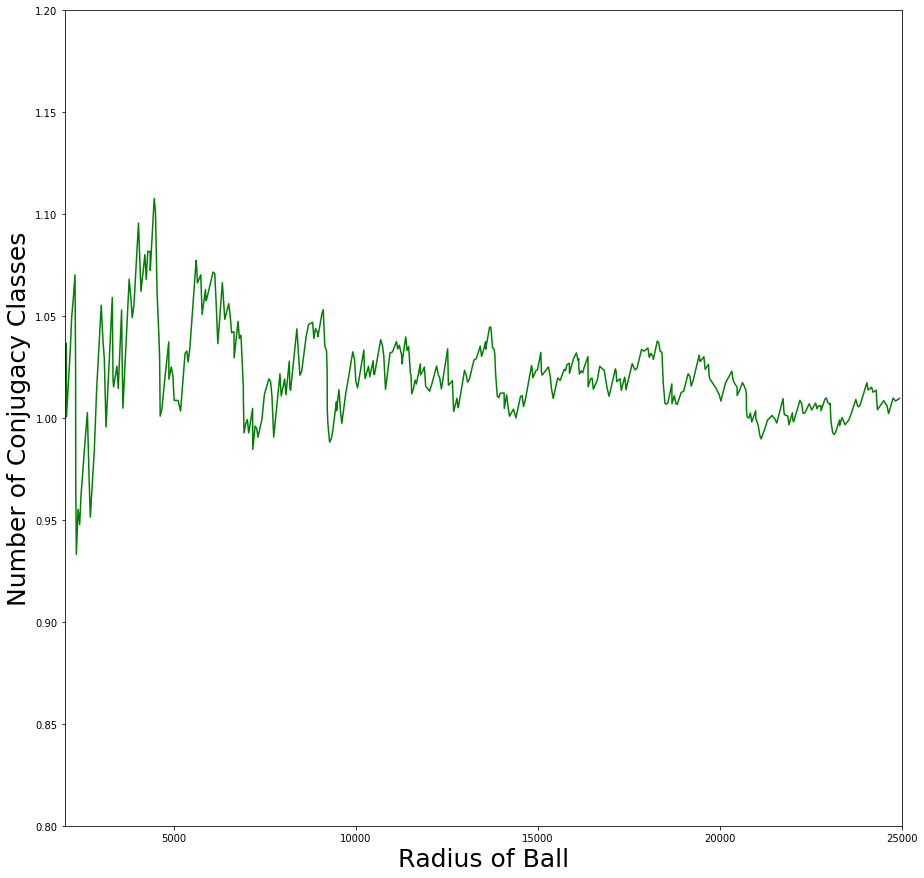}
\caption{A plot of $(T, \sum_{|t|\leq T} h'(t))$ for $T<25000$ vs the ratio of Epstein's estimate and $ \sum_{|t|\leq T} h'(t)$}
\label{ep}
\end{figure}

When we are done we get that $\mathcal{A}^{m_1} \mathcal{B}^{n_1}\cdots \mathcal{A}^{m_k} \mathcal{B}^{n_k}g=I$ and so $g=(\mathcal{A}^{m_1} \mathcal{B}^{n_1}\cdots \mathcal{A}^{m_k} \mathcal{B}^{n_k})^{-1}.$ Summing the exponents of the $\mathcal{A}$'s and $\mathcal{B}$'s will tell us if $g$ is in the commutator subgroup. And now that we actually know how to calculate $h'(t)$, we can look at the growth rate and compare it to Theorem \ref{growth} as seen in Figure \ref{ep}.


\section{Studying Conjugacy classes that lie in $L$}\label{ch:walks}

Recall that for $g\in H$ there exists $h\in H'$ such that $g$ can be rewritten as $g=\mathcal{A}^{m}\mathcal{B}^{n}h$ (see Lemma \ref{appendixlemma}). And the homology map, $\pi:H\rightarrow \mathbb{Z}^2$, is $\pi(g)=\left(m, n\right).$ So when $g\in H',$ $\pi(g)=(0,0)$ and this is the trivial homology.

For $g\in H,$ we represent its $\SL_2(\mathbb{Z})$ conjugacy class as $[g]=\{hgh^{-1} \mid h\in \SL_2(\mathbb{Z})\}$ and its $H$ conjugacy class as $[g]_H=\{hgh^{-1} \mid h\in H\}.$ While the homology map is clearly the same for all elements in a given $H$ conjugacy class, it is not true with respect to the $\SL_2(\mathbb{Z})$ conjugacy class. Meanwhile, word length is not even consistent with respect to the $H$ conjugacy class.

Recall that any $g\in H$ can also be written as $g=\mathcal{B}^n \delta \mathcal{A}^m$ for some $\delta=\mathcal{A}^{m_1}\cdots B^{n_k}\in \Delta.$ Then in Section \ref{ch:Pres}, we defined the word length of $g$ as $\ell(g)=\sum |m_i|+\sum |n_i|+|n|+|m|$ and the narrow length of $\delta$ as $\ell_N(\delta)=k$. 

\begin{example}\label{ex1walks}
 Consider $g_1=\left(\mathcal{A}\mathcal{B}^{-1}\right)^{10}\mathcal{A}^{-2}\mathcal{B}^2.$  For $h=\begin{pmatrix}1 & 1 \\0& 1 \end{pmatrix}\not\in H$, by definition, $hg_1h^{-1}\in [g_1]$ and one can compute $hg_1h^{-1}=\mathcal{A}\mathcal{B}^{10}\mathcal{A}^{-2} \mathcal{B}^{-1}{A} \mathcal{B}^{-1}.$ Here $\ell(g_1)=24$ while $\ell(hg_1h^{-1})=17$ and, similarly, $\pi(g_1)=(8,-8)$ while $\pi(hg_1h^{-1})= (0,8).$
\end{example}

Note, also, that the homology map is consistent on $\SL_2(\mathbb{Z})$ conjugacy classes when restricted to the trivial homology, but again the word length is not.

\begin{example}
 Now look at $g_2=\left(\mathcal{A}\mathcal{B}^{-1}\right)^{10}(\mathcal{A}^{-1}\mathcal{B})^{10}.$ Then $hg_2h^{-1}=\mathcal{A}\mathcal{B}^{10}\mathcal{A}^{-1}\mathcal{B}^{-10}$ for $h$ as in Example \ref{ex1walks}. Here  $\ell(g_2)=40$ while $\ell(hg_2h^{-1})=22$ and $\pi(g_2)=(0,0)=\pi(hg_2h^{-1}).$
\end{example}

\subsection{Properties of Length and Homology classes}

Since we have been looking at hyperbolic conjugacy classes with respect to $\SL_2(\mathbb{Z})$ that lie in $L$, we need to understand more about $\{\ell(hgh^{-1}) \mid h\in \SL_2(\mathbb{Z})\} \text{ and }\{\pi(hgh^{-1}) \mid h\in \SL_2(\mathbb{Z})\}$ for $g\in L.$ To do this, it is important to record some properties of the word length and the homology class of a given $g.$

\begin{lem}\label{somesmallproperties}
For $g\in H$, the following are true:
\begin{enumerate}
    \item $[g^t]=[g^{-1}]$
    \item $\ell(g)=\ell(g^t)=\ell(g^{-1})$
    \item $\pi(g)=\pi(g^t)^\tau = -\pi(g^{-1})$ (where $(a,b)^\tau=(b,a)$)
\end{enumerate}
\end{lem}
\begin{proof}
(1) is just due to the fact that $$\begin{pmatrix}0& -1\\1 & 0\end{pmatrix}\begin{pmatrix}a & c \\b & d \end{pmatrix}\begin{pmatrix}0& 1\\-1 & 0\end{pmatrix}=\begin{pmatrix} d & -b \\ -c & a \end{pmatrix}.$$

For (2) and (3) we write $g$ as $g=\mathcal{B}^{n_{k+1}}\mathcal{A}^{m_1}\mathcal{B}^{n_1}\cdots \mathcal{A}^{m_k}\mathcal{B}^{n_k}\mathcal{A}^{m_{k+1}}$ for $m_i,n_i\in \mathbb{Z}_{\neq 0}$ and $m_{k+1},n_{k+1}\in \mathbb{Z}.$ Then, we can rewrite $g^t$ and $g^{-1}$ as follows:
\begin{equation*}
    g^t=\mathcal{B}^{m_{k+1}}\mathcal{A}^{n_k}\mathcal{B}^{m_k}\cdots \mathcal{A}^{n_1}\mathcal{B}^{m_1}\mathcal{A}^{n_{k+1}} \text{ and } g^{-1}=\mathcal{A}^{-m_{k+1}}\mathcal{B}^{-n_k}\mathcal{A}^{-m_k}\cdots \mathcal{B}^{-n_1}\mathcal{A}^{-m_1}\mathcal{B}^{-n_{k+1}}.
\end{equation*} 

Thus it is clear that since, $\ell(g)=\sum |m_i|+ \sum |n_i|$, so does $\ell(g^{-1})$ and $\ell(g^t).$ It is similarly easy to see that $(\sum m_i, \sum n_i)=(\sum n_i, \sum m_i)^\tau=-(\sum -m_i, \sum -n_i)$ and so we are done.  
\end{proof}

\begin{lem}
Let $g\in \Delta$, $g\neq I,$ and $h\in H.$ If $hgh^{-1}\in \Delta$ then $h\in \Delta \cup \Delta^{-1}$ and $hgh^{-1}$ is a cyclic rotation of $g.$
\end{lem}
For the proof of this lemma, see pages 80-81 in \cite{dissertation}. It uses a proof by contradiction for $h\in \Delta\cup \Delta^{-1}$ and then inductive argument on the narrow length of $h$ to prove the cyclic rotation of $g$. Now by definition of the word length and narrow length we can conclude the following.
\begin{cor}\label{cyclicrotation}
Let $g\in \Delta$ and $h\in H.$ If $hgh^{-1}\in \Delta$ then $\ell(g)=\ell(hgh^{-1})$ and $\ell_N(g)=\ell_N(hgh^{-1}).$
\end{cor}

\begin{lem}\label{randoquick}
If $g\in L,$ then $\pi(g)\equiv (0,0) \bmod{4}.$
\end{lem}
\begin{proof}
Since $g\in H, $  we can write $g=\mathcal{A}^m\mathcal{B}^n h$ for some $h\in H'.$ So $\pi(g)=(m,n)$ and all that remains is to show that $4|m$ and $4|n.$ Given that $g\in L$ and $h\in L,$ we get that
\begin{equation*}
 \mathcal{A}^m\mathcal{B}^n=\begin{pmatrix}1 & 2m \\0 & 1 \end{pmatrix} \begin{pmatrix}1 & 0 \\2n & 1 \end{pmatrix} = \begin{pmatrix}1+4m n & 2m \\ 2n & 1 \end{pmatrix}\in L.
\end{equation*} In order for it to be in $L$, however, $8|2m$ and $8|2n.$ Thus we are done.
\end{proof}

We will now look at how $\SL_2(\mathbb{Z})/\Gamma(2)$ acts on $L$. We use the following representation of $\SL_2(\mathbb{Z})/\Gamma(2)$: \begin{equation}\label{star}
    \left\{\begin{pmatrix}1 & 0\\0 & 1\end{pmatrix}, \begin{pmatrix}1 & 1\\0 & 1\end{pmatrix}, \begin{pmatrix}1 & 0\\1 & 1\end{pmatrix},\begin{pmatrix}1 & -1\\1 & 0\end{pmatrix},\begin{pmatrix}0 & -1\\1 & 1\end{pmatrix}, \begin{pmatrix}0 & -1\\1 & 0\end{pmatrix}\right\}  =\{\gamma_0, \gamma_1, \gamma_2, \gamma_3, \gamma_4, \gamma_5\}.
\end{equation}

\begin{lem}\label{finalposlem}
Let $g\in L$ and $\pi(g)=(m,n).$ Then, for $\gamma_i$ defined in Equation \ref{star},
\begin{enumerate}
    \item $\pi(\gamma_1 g\gamma_1^{-1})=(m+n,-n)$,
    \item $\pi(\gamma_2 g\gamma_2^{-1})=(-m,m+n)$,
    \item $\pi(\gamma_3 g\gamma_3^{-1})=(-(m+n),m)$, 
    \item $\pi(\gamma_4 g\gamma_4^{-1})=(n,-(m+n))$,  
    \item $\pi(\gamma_5 g\gamma_5^{-1})=(-n,-m)$, and 
    \item $\pi(g^t)=(n,m).$
\end{enumerate}
\end{lem} 
The proof of which is a direct computation for $g=\mathcal{A}^{m_1}\mathcal{B}^{n_1}\cdots \mathcal{A}^{m_k}\mathcal{B}^{n_k}. $ For example conjugating $\mathcal{A}$ by $\gamma_1$ is $\mathcal{A}$  and  $\gamma_1 \mathcal{B} \gamma_1^{-1}= -\mathcal{A}\mathcal{B}^{-1}.$
Keeping in mind that $\sum n_i$ is even (Lemma \ref{randoquick}) and so the negatives can be ignored, we get
\begin{equation}
\pi \left(\gamma_1 g\gamma_1^{-1}\right)=\pi \left( \mathcal{A}^{m_1}\left(\mathcal{A}\mathcal{B}^{-1}\right)^{n_1}\cdots \mathcal{A}^{m_k}\left(\mathcal{A}\mathcal{B}^{-1}\right)^{n_k}    \right)= \left(\sum m_i + \sum n_i , -\sum n_i\right) = (m+n, -n).
\end{equation}

 We can take the representatives of the $\SL_2(\mathbb{Z}$) conjugacy classes we found in Section \ref{ch:existence} and get the $H$ conjugacy classes by conjugating each matrix by the $\gamma_i$'s in Equation \ref{star} and then comparing their decompositions into $\mathcal{A}$'s and $\mathcal{B}$'s. Since $\pi(h)$ is the same for all $h\in [g]_H,$ we let $\pi([g]_H):=\pi(g)$ and we look at the heat maps the set of $H$ conjugacy classes produce.

 \begin{cor}\label{mm}
Suppose $g\in L$ and that $g^t\in [g]_H.$ Then $\pi(g)=(m,m)$ for some $m\in \mathbb{Z}.$
\end{cor}

 We can take the representatives of the $\SL_2(\mathbb{Z}$) conjugacy classes we found in Section \ref{ch:existence} and get the $H$ conjugacy classes by conjugating each matrix by the $\gamma_i$'s in Equation \ref{star} and then comparing their decompositions into $\mathcal{A}$'s and $\mathcal{B}$'s. Since $\pi(h)$ is the same for all $h\in [g]_H,$ we let $\pi([g]_H):=\pi(g)$ and we look at the heat maps the set of $H$ conjugacy classes produce.

\begin{figure}[h!]
\center
\includegraphics[width=0.45\textwidth]{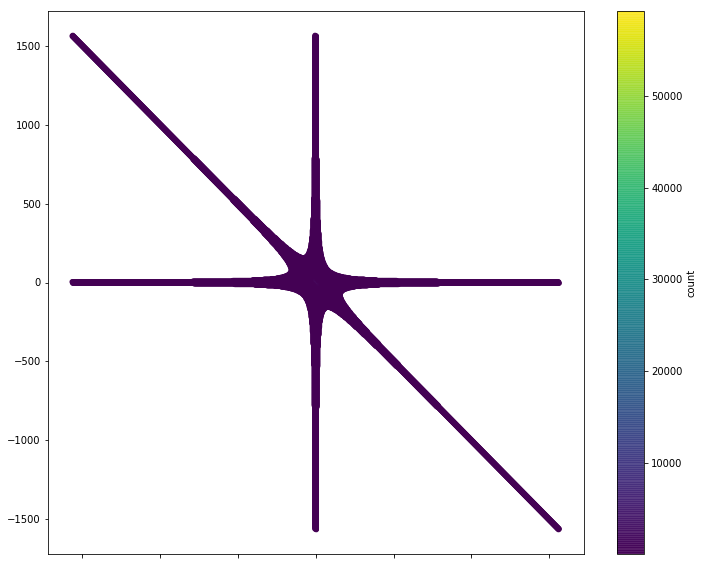}
\includegraphics[width=0.45\textwidth]{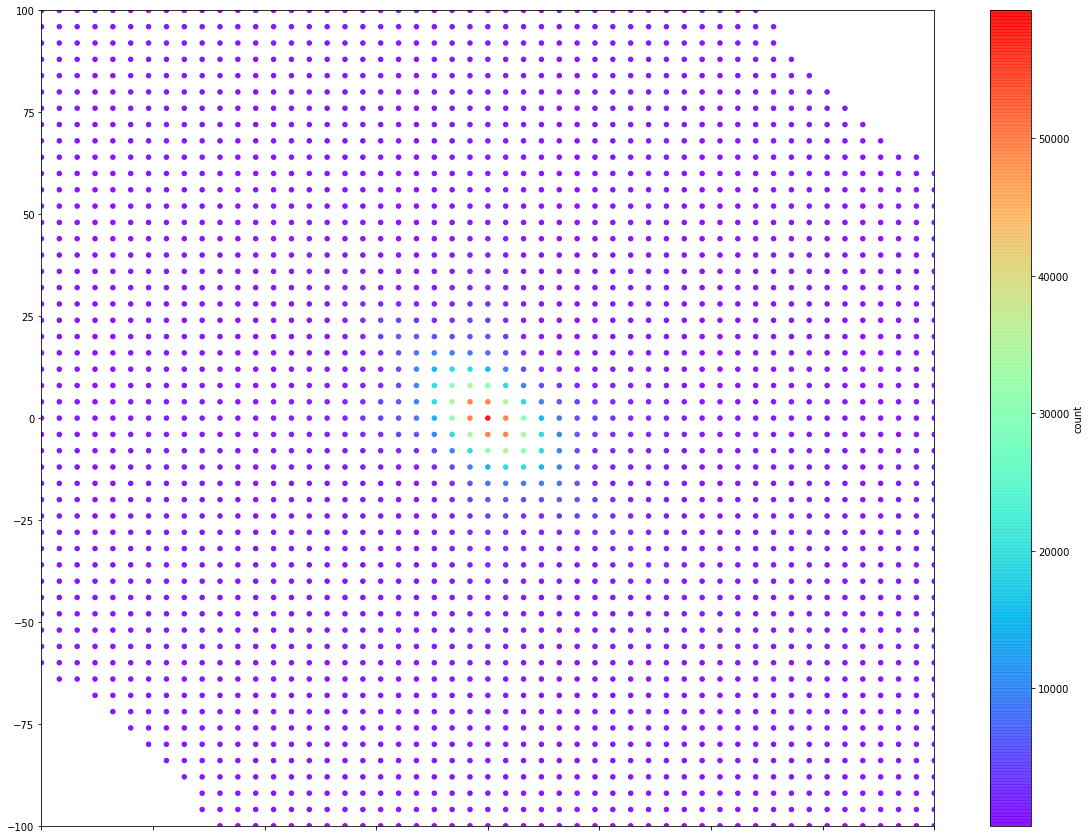}
\caption{A heat map of $\pi([g]_H)$ for $g\in L$ with $|\Tr(g)|<25000$ zoomed in.}
\label{multofwalks}
\end{figure}

The heat map of $\pi([g]_H)$ is symmetric about the line $y=x$ and about the origin. Let $g$ be a representative of an $H$ conjugacy class. Either $g^t \in [g]_H,$ which by Corollary \ref{mm} implies that $\pi([g]_H)=(m,m),$ or it is not. In the latter case, $\pi([g^t]_H)=\pi(g^t)=\pi(g)^\tau.$ Thus for any $(m,n)$ plotted, $(n,m)$ is also and so the map is symmetric about $y=x$. Similarly, if $g^{-1}\in [g]_H,$ then $\pi(g)=\pi(g^{-1})$ and so $\pi([g]_H)=(0,0).$ If $g^{-1}\not\in [g]_H,$ then $\pi([g^{-1}]_H)=(-m,-n)$ is plotted and so the heatmap is symmetric about the origin.

\begin{figure}[h!]
\center
\includegraphics[width=0.45\textwidth]{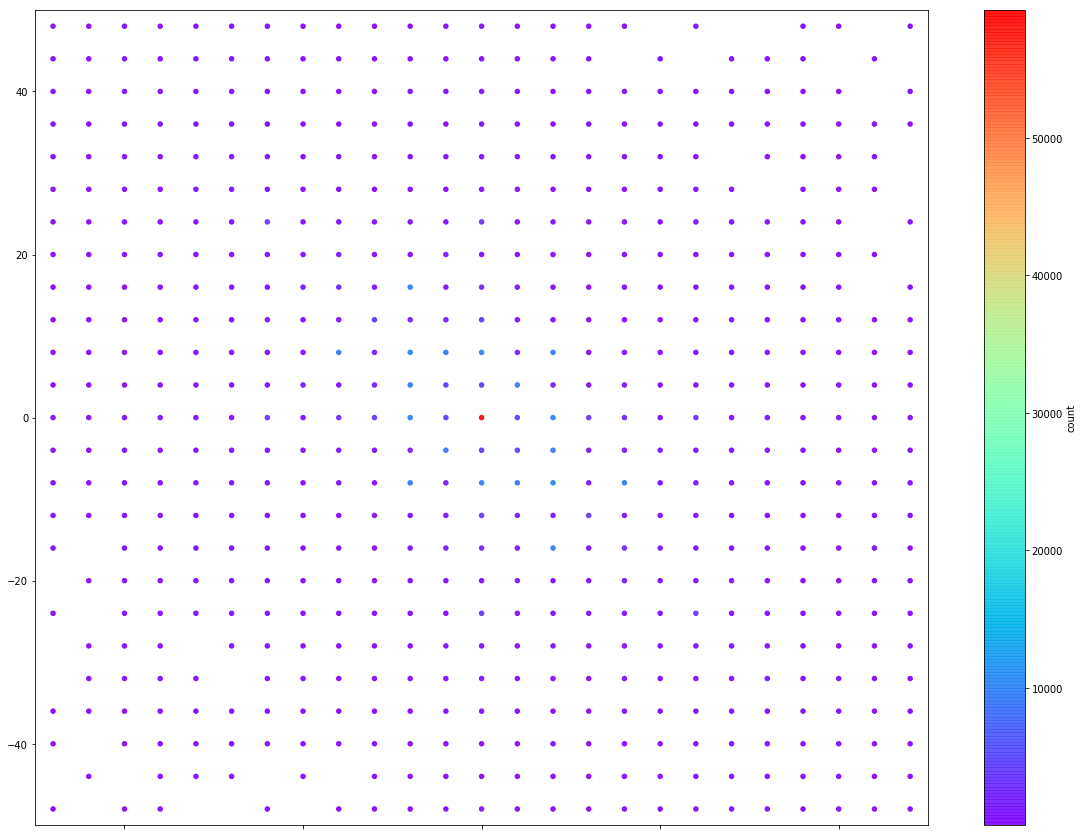}
\includegraphics[width=0.45\textwidth]{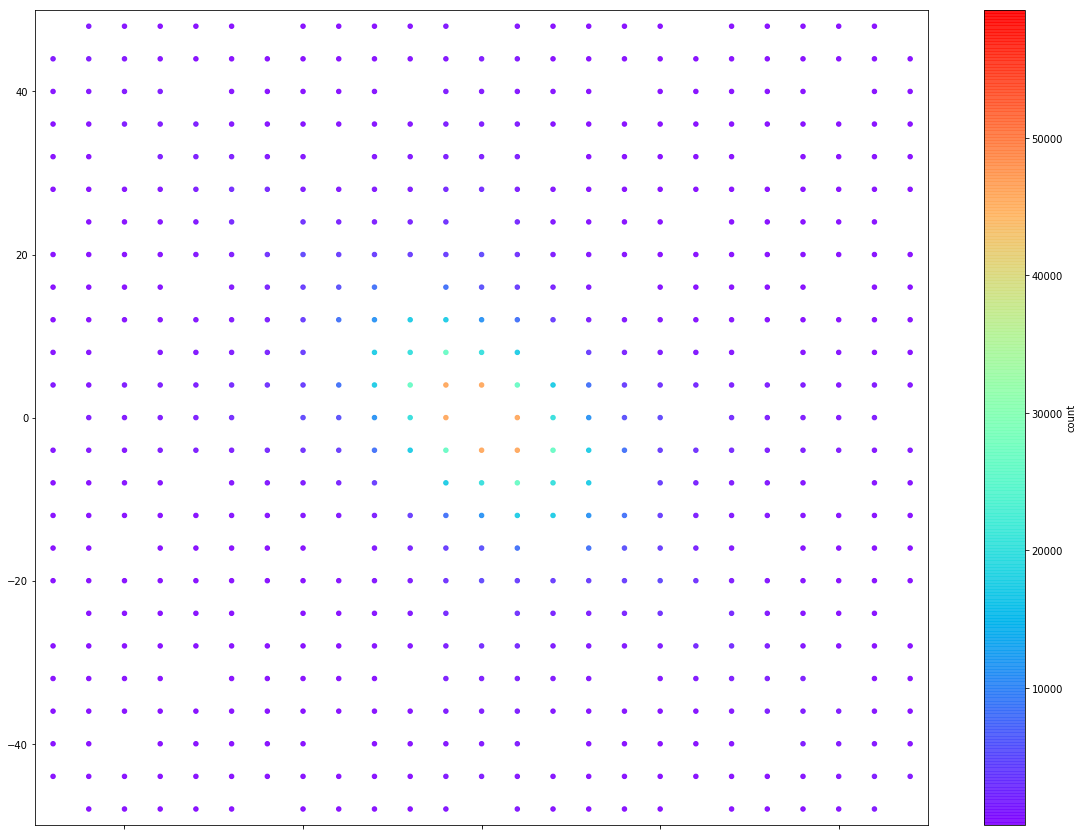}
\caption{A heat map of $\pi([g]_H)$ for $g\in L$ with $|\Tr(g)|<25000.$ On the left is admissible traces and on the right is nonadmissible traces.}
\label{allwalksBig}
\end{figure}

\begin{lem}\label{8m8m24n}
For any $g\in H$, if there exists $m,n\in \mathbb{Z}$ such that $\pi(g)=(8m, 8m+24n)$ then $\Tr(g)$ is admissible.
\end{lem}
\begin{proof}
Take any $m,n\in \mathbb{Z}$ then write $g$ as $g=\mathcal{A}^{8m}\mathcal{B}^{8(3n+m)} h$ where $h\in H'.$ Then we just need to show that $\Tr(g)$ is admissible. It is fairly easy to see (by looking at $\Gamma(2)' \bmod{16}$) that there are two forms that the matrices can take, \begin{equation*}
    \begin{pmatrix} 1+8x & 16y \\ 16z & 1+64v-8x
  \end{pmatrix} \textit{ and }
    \begin{pmatrix} 5+8x & 16y+8 \\ 16z+8 & 13+64v-8x
  \end{pmatrix}.
  \end{equation*} The first form is when traces are congruent to 2 or 66 modulo 256 and the second is when the trace is congruent to 18 or 146 modulo 256.

If $h$ is of the first form we see that,
\begin{equation*}
 \Tr\left(\mathcal{A}^{8m}\mathcal{B}^{8(3n+m)}h\right)= \Tr\left(\begin{pmatrix}1+256m(3n+m) & 16m \\ 16(3n+m) & 1 \end{pmatrix}\begin{pmatrix} 1+8x & 16y \\ 16z & 1+64v-8x
  \end{pmatrix}\right) \equiv \Tr(h) \bmod 256.
 \end{equation*}
Else, $h$ is of the second type and 
 \begin{align*}
 \Tr\left(\mathcal{A}^{8m}\mathcal{B}^{8(3n+m)}h\right)&= \Tr\left(\begin{pmatrix}1+256m(3n+m) & 16m \\ 16(3n+m) & 1 \end{pmatrix}\begin{pmatrix} 5+8x & 16y+8 \\ 16z+8 & 13+64v-8x
  \end{pmatrix}\right)\\&\equiv \Tr(h)+128s \bmod 256
 \end{align*} for some appropriate $s$. So multiplying $\mathcal{A}^{8m}\mathcal{B}^{8(3n+m)}$ to either of the two forms shows that the trace stays in $\{2, 18, 66, 146\} \bmod{256}$. 
 
 Now we need to do the same argument modulo 9 which is not as clean. We want to show that for any $h\in \Gamma(2)' \bmod{9}$, $\Tr(\mathcal{A}^{8m}\mathcal{B}^{8(3n+m)}h) \bmod{9}\in \{0,2,3,6,7\}.$ For this, we can just do an exhaustive search. First we take the 216 matrices in $\Gamma(2)' \bmod{9}$ (see Lemma \ref{big}), then multiply them on the left by $\mathcal{A}^{8m}\mathcal{B}^{8(3n+m)}.$ For this list of matrices, we evaluate the trace modulo 9 for all $0\leq n,m<9$ and see that the set of possible traces is just $\{0,2,3,6,7\}$ thus giving us that the trace must be admissible, though, not necessarily the same trace as $h.$
\end{proof}

\subsection{Choosing ``Good" Representatives}

Now that we know these properties of word length and homology we will extend them to conjugacy classes. 
\begin{defn}
For $g\in L$, the word length of its $\SL_2(\mathbb{Z})$ conjugacy class and its $H$ conjugacy class are defined as
$$\ell([g])=\min_{h\in [g]} \ell(h) \textit{ and } \ell([g]_H)=\min_{h\in [g]_H} \ell(h).$$
\end{defn}

Notice here that for any $h\in [g]_H\cap \Delta,$ $\ell([g]_H)=\ell(h)$ via Corollary \ref{cyclicrotation} and also, by definition, $\ell([g])\leq \ell([g]_H).$

\begin{lem}
If $\ell(g)=\ell([g]),$ then $\ell(g^t)=\ell([g^t]).$
\end{lem}
\begin{proof}
Suppose not, that is suppose there exists $h\in [g^t]$ such that $\ell(h) < \ell(g^t).$ Then, since $h=\gamma g^t \gamma^{-1}$ for some $\gamma\in \SL_2(\mathbb{Z}),$ $h^t=(\gamma^{-1})^t g \gamma^t$ which implies that $h^t\in [g],$ but then we get a contradiction since, by Lemma \ref{somesmallproperties}, $\ell(g)\leq \ell(h^t)=\ell(h)<\ell(g^t).$
\end{proof}

The idea is now to modify the representatives of the conjugacy classes (produced by the algorithm described in Section \ref{algoExistence}) to force the following three conditions on the representatives chosen:
\begin{enumerate}
    \item $g\in \Delta \cap L$
    \item $\ell(g)=\ell([g]),$
    \item If $g$ is a representative of the conjugacy class and $g^t\not\in [g],$ then $g^t$ is the chosen representative for $[g^t]$.
\end{enumerate} 
We then define $\Tilde{\rho}([g])=g$ where $g$ is this ``good" representative (in the sense that we get from the modified list of representatives). This representation is chosen so that the word length is as short as possible. Later on in this section, and the next, we will observe that having a small word length is helpful and grouping elements in terms of their $\SL_2(\mathbb{Z})$ conjugacy classes accomplishes this.

Let $g\in L$ and $\Tr(g)\neq 2.$  Then, we can define a mapping, $\sigma:\textit{hyperbolic matrices in } L \rightarrow \Delta,$ as follows. Let $\Tilde{h}\in H$ be such that $\Tilde{h}g\Tilde{h}^{-1}$ is the first counterclockwise cyclic rotation of $g$ that lies in $\Delta,$ then $\sigma(g)=\Tilde{h}g\Tilde{h}^{-1}.$

\begin{example}
Let $g=\mathcal{A}^2\mathcal{B}^4\mathcal{A}^2\mathcal{B}^{-4}\mathcal{A}^4$, then $\sigma(g)=\mathcal{A}^2\mathcal{B}^{-4}\mathcal{A}^6\mathcal{B}^4$ after conjugation by the element $\Tilde{h}=\mathcal{B}^{-4}\mathcal{A}^{-2}.$ Note that if we had instead conjugated by $\mathcal{A}^4$ we would have gotten $\mathcal{A}^6\mathcal{B}^4\mathcal{A}^2\mathcal{B}^{-4}.$ We don't choose to do this even though $\mathcal{A}^{4}$ has smaller word length than $\Tilde{h}$, since we would be rotating $g$ clockwise.
\end{example}

\begin{lem}\label{ohlala} If $g\in \Delta\cap L$ then for some $\gamma_i \in \SL_2(\mathbb{Z})/\Gamma(2),$ $\ell(\sigma(\gamma_i g \gamma_i^{-1}))=\ell([g]).$
\end{lem}
\begin{proof}
By definition of $\ell([g])$ and what it means to be in the conjugacy class, there exists an element, $k\in \SL_2(\mathbb{Z}),$ such that $kgk^{-1}\in \Delta \cap L$ and $\ell(kgk^{-1})=\ell([g]).$ By definition of coset representations, there exists $\gamma \in \SL_2(\mathbb{Z})/\Gamma(2)$ and $h\in \Gamma(2)$ such that $k=h\gamma.$  Thus, $\ell(h\gamma g \gamma^{-1} h^{-1})=\ell([g]).$ Similarly, $\sigma(\gamma g \gamma^{-1})\in \Delta \cap L.$ Note that $\sigma(\gamma g \gamma^{-1})$ and $kgk^{-1}$ are both in $\Delta$ and, if we write $\sigma(\gamma g \gamma^{-1})=\Tilde{h}\gamma g \gamma^{-1}\Tilde{h}^{-1}$, then $(h\Tilde{h}^{-1})\sigma(\gamma g \gamma^{-1})(h\Tilde{h}^{-1})^{-1}=kgk^{-1}.$ So Corollary \ref{cyclicrotation} gives us that $\ell(\sigma(\gamma g\gamma^{-1}))=\ell([g]).$
\end{proof}

Now we are ready to describe how we choose our final representatives of the hyperbolic conjugacy classes of $\SL_2(\mathbb{Z})$ that lie in $L$. The goal is to find representatives of the conjugacy classes with the shortest word length, which, by Lemma \ref{ohlala}, means that we only need to check $\sigma(\gamma_i g\gamma_i^{-1})$ for $\gamma_i$ as in Equation \ref{star}. However $\ell(\sigma(\gamma_5 g\gamma_5^{-1}))=\ell(g)$ via direct caluclations and so we do not need to check this case.

\begin{figure}[h!]
\center
\includegraphics[width=0.45\textwidth]{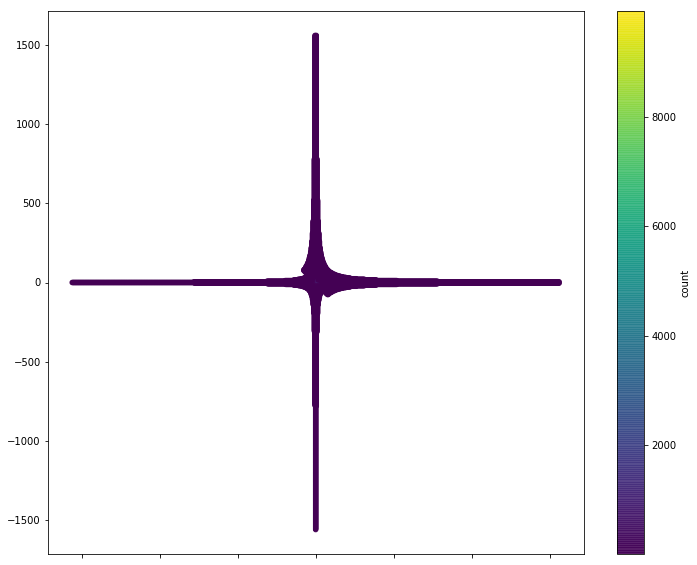}
\includegraphics[width=0.45\textwidth]{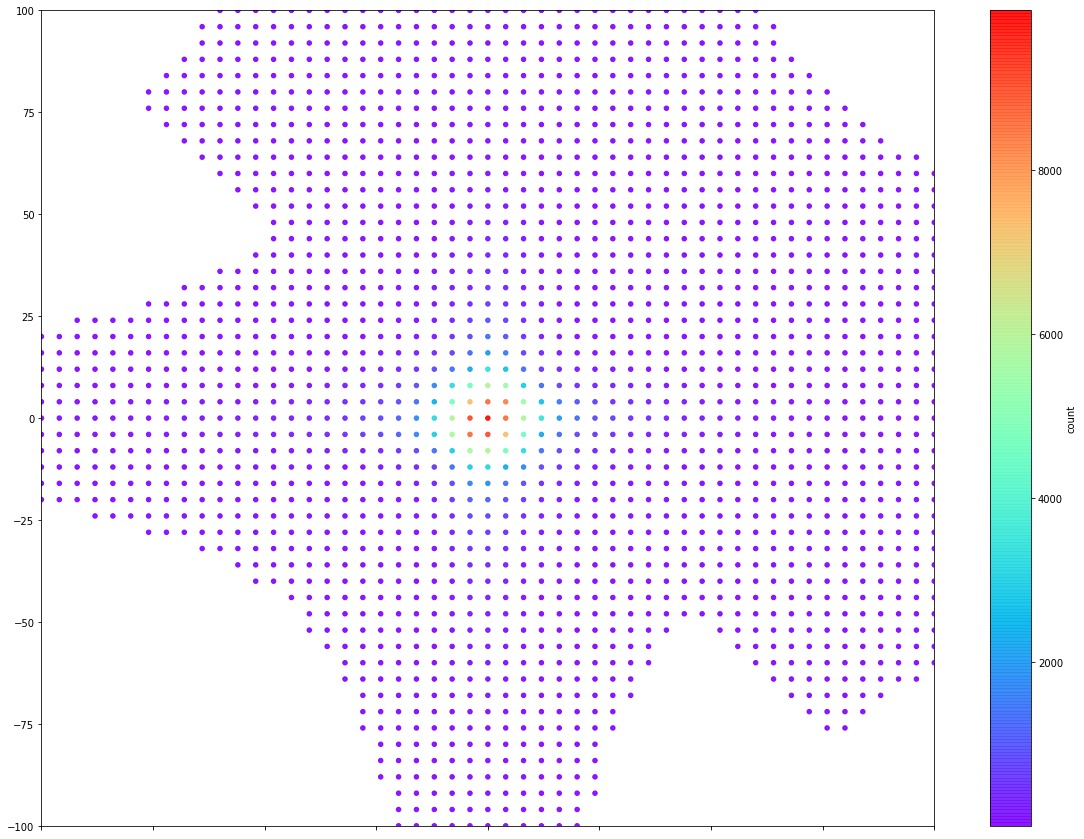}
\caption{A heat map of $\pi(\Tilde{\rho}([g]))$ with $|\Tr(g)|<25000$.}
\end{figure}

We take the representatives gained from Section \ref{algoExistence} with trace $t$ and name this set $\hat{h}(t)$. For each representative, $g\in \hat{h}(t)$, we first let $h=\sigma(\gamma_0 g\gamma_0^{-1})=\sigma(g)$ and then we compute $\sigma(\gamma_i g\gamma_i^{-1})$ for $i$ ranging from 1 to 4. If $\ell(\sigma(\gamma_i g\gamma_i^{-1}))<\ell(h)$, we replace $h$ by $\sigma(\gamma_i g\gamma_i^{-1}).$ Once completed, we replace $g$ with $h$ in $\hat{h}(t)$. Next, we check if $g^t\in [g]$. If it is not, we find the representative for $[g^t]$ in $\hat{h}(t)$ and replace it with $h^t$.

\subsection{A Lower Bound on the Length $g$}

\begin{figure}[H]
\center
\includegraphics[width=0.4\textwidth]{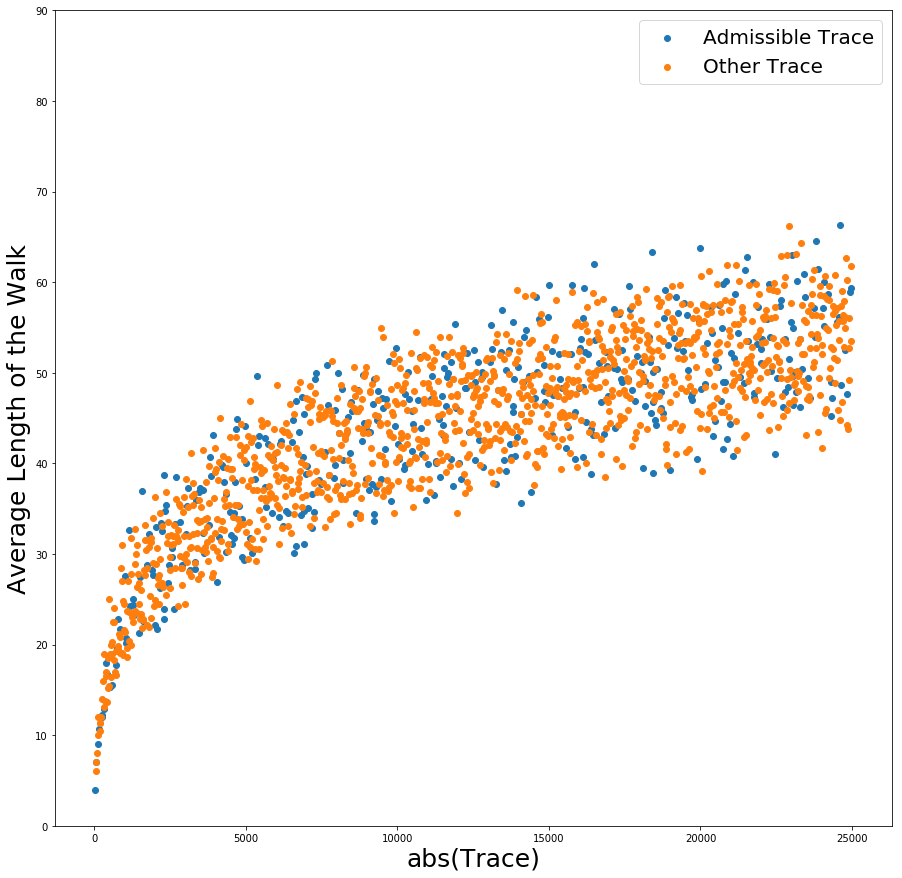}
\includegraphics[width=0.4\textwidth]{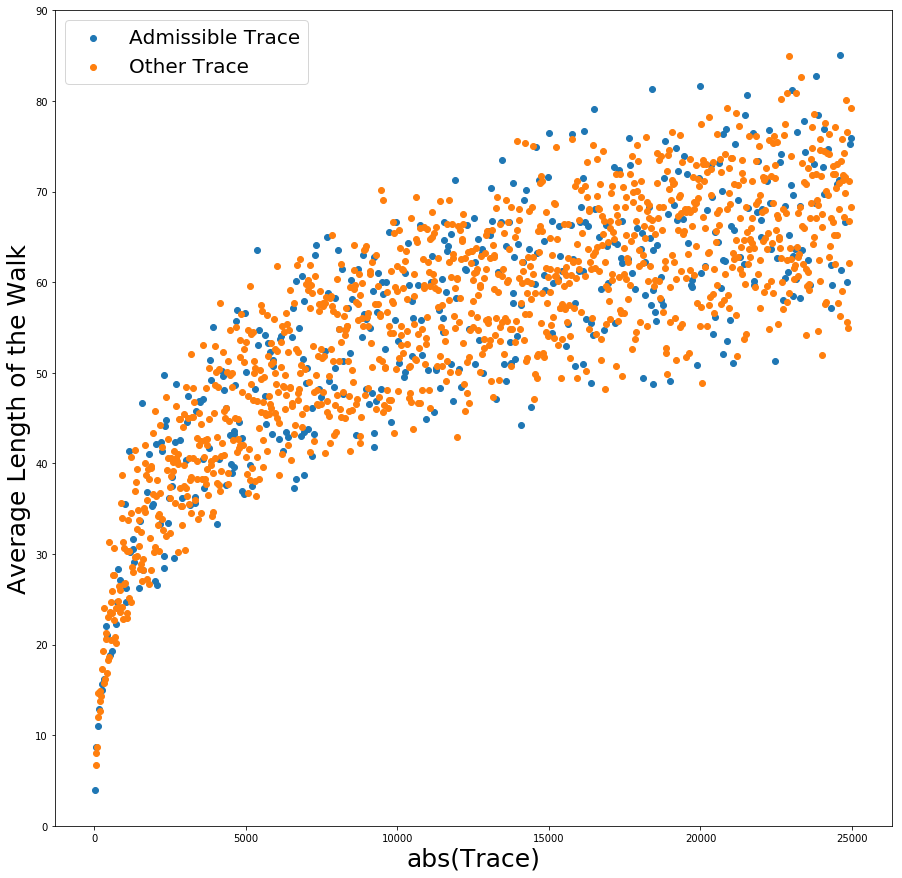}
\caption{Average of $\ell([g])$ over all hyperbolic conjugacy classes in $L$ with trace $t$ for $|t|<25000$ on the left and that of $\ell([g]_H)$ on the right.}
\label{htpht}
\end{figure}

\begin{conj}\label{iwish}
As $t$ (admissible) approaches infinity, $$\frac{1}{h(t)}\sum_{\substack{[g]\subset L\\ \Tr(g)=t}} \ell([g])\asymp \log(|t|).$$
\end{conj}
Though we cannot prove this, we can (and do in pages 82-84 of \cite{dissertation}) prove the following theorem via induction after first proving that  $\log\left(\Tr\left(\begin{pmatrix}5 &2\\2&1\end{pmatrix}^k\right)-2\right)\leq 2k.$

\begin{thm}
If $g\in \Delta$ and $g\neq I$, then $\ell(g)>\log(|\Tr(g)|/2)$.
\end{thm} 
This theorem gives the lower bound for Conjecture \ref{iwish} but there is no analogy for the upper bound given that the word length can get as large as the trace. The goal, then, is to prove that these don't happen often, which based on the figure in this section, seems like a reasonable conjecture. Truly though, we want to know how the word length of a conjugacy class behaves on average, since these lengths can get quite large with respect to $t$. In the following we can see one example of just how big it can get. 
\begin{example}
Take $t=-92142$ which is admissible. Then there exists a matrix $g_1\in L$ such that $\Tr(g_1)=t$ and $\ell(g_1)=26$. This $g_1$ is the following:
\begin{equation*}
  g_1=\begin{pmatrix}-77651& 30920\\36392& -14491 \end{pmatrix}=\mathcal{A}^{-1}\mathcal{B}^{-4}\mathcal{A}\mathcal{B}^{-6}\mathcal{A}\mathcal{B}^{11}\mathcal{A}^{-1}\mathcal{B}^{-1}.
\end{equation*} Note, $g_1\in H'.$ 
However, in a different conjugacy class, there exists $g_2\in L$ that has $\Tr(g_2)=t$ and $\ell(g_2)=5762$. It is
\begin{equation*}
 g_2=\begin{pmatrix}-115179 & 46072\\-57592& 23037 \end{pmatrix}=  \mathcal{A}\mathcal{B}^{-5759}\mathcal{A}^{-1}\mathcal{B}^{-1}
\end{equation*} and clearly $g_2\not\in H'.$
\end{example}

This example demonstrates that there exist traces where the word length can get very large but hopefully, when we limit the search to elements in the commutator subgroup, they do not grow as quickly. 

\subsection{Looking at a random walk}

\begin{figure}[h]
\center
\includegraphics[width=0.35\textwidth]{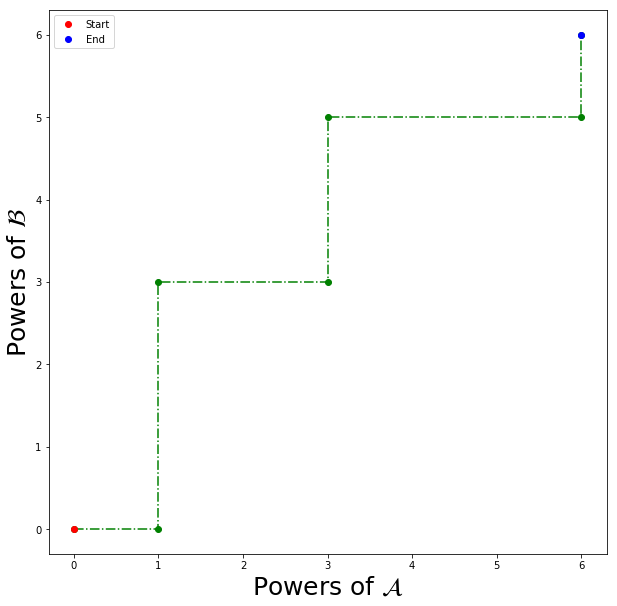}
\includegraphics[width=0.35\textwidth]{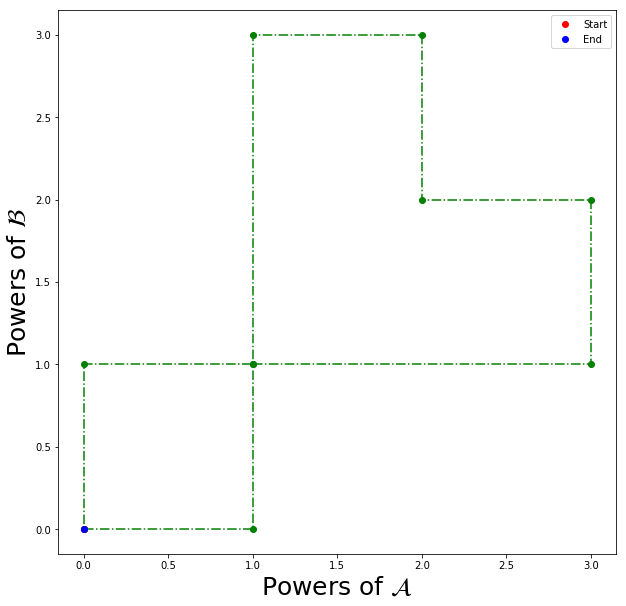}
\caption{The walk on the left is for the matrix $\mathcal{A}^{1}\mathcal{B}^{3}\mathcal{A}^{2}\mathcal{B}^{2}\mathcal{A}^{3}\mathcal{B}^{1}.$ The walk on the right hand side is for the matrix $\mathcal{A}\mathcal{B}\mathcal{A}^{2}\mathcal{B}\mathcal{A}^{-1}\mathcal{B}\mathcal{A}^{-1}\mathcal{B}^{-2}\mathcal{A}^{-1}\mathcal{B}^{-1}.$}
\label{showingawalk}
\end{figure}

For $g\in \Delta$ where $g=\mathcal{A}^{m_1}\mathcal{B}^{n_1}\cdots \mathcal{A}^{m_k}\mathcal{B}^{n_k}$, a walk is the set of points, $$
\left\{(0,0), (m_1, 0), (m_1, n_1), (m_1+m_2, n_1), \ldots, \left(\sum m_i, \sum n_i\right)\right\}.$$ Then the walk ends at $\pi(g)$ and the length of the walk is $\ell(g)$ (see Figure \ref{showingawalk}). It is for this reason that we can model our set of hyperbolic conjugacy classes as non-backtracking random walks in the two dimensional lattice and look for when walks end at the origin (points such that $\pi(g)=(0,0)$ and are thus in $H'$).

Now we will look at the probability of a non-backtracking random walk in 2 dimensions ending at the origin. We noted that the $h(t)$ representatives of hyperbolic conjugacy classes appear to have an average word length on the order of $\log(|t|)$. It is known that the probability of a random walk on a 2d lattice of length $2N$ ending at the origin is asymptotically $(\pi N)^{-1}.$ This calculation can be done by first noting that the total number of walks is $4^{2N}$ and the number of walks that end at the origin is $\sum \binom{2N}{2k}\binom{2k}{k}\binom{2N-2k}{N-k}.$ The latter value comes from noticing that, in order to return to the origin, the total amount of times one travels up is equal to the total number of times one travels down and the total number of times one travels left must equal the times traveling right. Using identities of binomial coefficients this probability reduces to $\frac{\binom{2N}{N}^2}{4^{2N}}$ and then Stirling's Approximation gives the asymptotic $(\pi N)^{-1}.$  The probability of a non-backtracking walk of length $2N$ ending at the origin is then also asymptotically $(c\pi N)^{-1}$ for some $c$ \cite{OEIS}. If the model accurately predicts the homology class, then for large enough admissible $t$, the size of $h'(t)$ should be of the order $h(t)/\log(|t|)$.




\section{Commutator Width}\label{ch:coms}

We might ask ourselves if $g\in \Gamma(2)'$ can it be written as a 1-commutator. We can further generalize this question to ask, what is the minimal number of commutators needed to write $g\in \Gamma(2)'.$ We call this value the commutator width of $g$.  Given that $\gamma[h,k]\gamma^{-1}=[\gamma h\gamma^{-1}, \gamma k \gamma^{-1}],$ it is clear that $g$ has the same width as any other element in its $\SL_2(\mathbb{Z})$ conjugacy class. So the question reduces, as usual, to answering the question in $\Delta\cap \Gamma(2)'$. 

For certain $g$, this question is not hard to answer as can be seen in Lemma \ref{easycase1} below.

\begin{lem}\label{easycase1}
Let $g\in \Delta\cap \Gamma(2)'$ and let $g$ have narrow length equal to 2 or 3. Then $g$ and $g^{-1}$ can both be written as a 1-commutator in $\Gamma(2)'.$
\end{lem}

\begin{proof}
Suppose $g$ has narrow length 2. Then, let $$g=\mathcal{A}^{m_1}\mathcal{B}^{n_1} \mathcal{A}^{m_2}\mathcal{B}^{n_2}=\mathcal{A}^{m_1}\mathcal{B}^{n_1}\mathcal{A}^{-m_1}\mathcal{B}^{-n_1}=[\mathcal{A}^{m_1},\mathcal{B}^{n_1}].$$ Suppose instead that $g$ has narrow length 3. Then, \begin{equation*}
  g=\mathcal{A}^{m_1}\mathcal{B}^{n_1}\mathcal{A}^{m_2}\mathcal{B}^{n_2}\mathcal{A}^{m_3}\mathcal{B}^{n_3}=[\mathcal{A}^{m_1}\mathcal{B}^{-n_2},\mathcal{B}^{n_1+n_2}\mathcal{A}^{m_2}].
\end{equation*}
Since $g$ must be a 1-commutator, $g=[X,Y]$ for some $X,Y\in H$ which then implies that $g^{-1}$ is also a 1-commutator given that $g^{-1}=[X,Y]^{-1}=[Y,X].$
\end{proof}
But the commutator widths are not always this easy to see.

\subsection{Algorithm for finding the width of an element in $\Gamma(2)'$}

The algorithm by Goldstein and Turner \cite{GoldsteinRichard1979Aotg} tells us not only if a conjugacy class can be written as a 1-commutator, but the minimum number of commutators we need to write it. The implementation of this algorithm can be found in the repository. \footnote{\href{https://github.com/BrookeOgrodnik/CommutatorSubgroups/blob/master/PythonCode/Genus\_of\_traces.ipynb}{https://github.com/BrookeOgrodnik/CommutatorSubgroups/blob/master/PythonCode/\newline Genus\_of\_traces.ipynb}}

\textit{A sketch of an Algorithm for finding the commutator width of an element in $\Delta \cap \Gamma(2)'$}
\begin{enumerate}
  \item Given $g\in \Delta \cap \Gamma(2)'$
  \item Write $g=\mathcal{A}^{m_1}\mathcal{B}^{n_1}\cdots \mathcal{A}^{m_k}\mathcal{B}^{n_k}$
  \item Convert into a list of $a$'s, $b$'s, $c$'s, $d$'s where $\mathcal{A}\rightarrow a$, $\mathcal{B}\rightarrow b$, $\mathcal{A}^{-1}\rightarrow c$, $\mathcal{B}^{-1}\rightarrow d$
  \item Pair the $a$'s and $c$'s together 
  \item Pair the $b$'s and $d$'s together 
  \item Assign labels to all of the pairings in order from left to right.
  \item Plot the points, in order they appear, on a circle and connect two points via a line segment if they were matched above.
  \item Create a matrix that at position $(i,j)$ is $1$ if the line associated with label $i$ intersects with line associated with label $j$ and $0$ else.
  \item Calculate the rank of the matrix in $\mathbb{Z}/2\mathbb{Z}$ and divide it by 2. 
  \item The minimum of these values is the width of $g$.
\end{enumerate}

For $g\in \Delta \cap \Gamma(2)',$ we can begin by searching for a mapping that has width 1. If we find one, then we can stop because we already hit the best case scenario. However, if $g$ has no width 1 mappings we must check every single mapping since we don't yet know if we can do better.  Given that $\ell(g)>\log(|t|/2)$, as the trace gets larger, this algorithm, obviously, becomes less helpful. 

\begin{remark}
Let $g\in \Delta \cap \Gamma(2)'$ and write it as $g=\mathcal{A}^{m_1}\mathcal{B}^{n_1}\cdots \mathcal{A}^{m_k}\mathcal{B}^{n_k}$ then at worst, the number of cases that would need to be checked is $\left(\frac{1}{2}\sum |m_i|\right)!\left(\frac{1}{2}\sum |n_i|\right)!$.
\end{remark}

Questions about commutator widths are common. As a comparison, it is known that the commutator subgroup of  $\SL_3(\mathbb{Z})$ is $\SL_3(\mathbb{Z})$ and this is not difficult to verify given \cite{MR1079696}. Also,  \cite{MR130917} gives us that $g\in \SL_3(\mathbb{Z}/p\mathbb{Z})'$ is a 1-commutator. It is then conjectured  that all $g\in \SL_3(\mathbb{Z})$ are 1-commutators. From previous theorems in Section \ref{ch:Admi} and, again, \cite{MR130917}, it is easy to see that, for our case, $g \in \Gamma(2)' \bmod{p}$ is a 1-commutator for prime $p\neq 3$ and from the examples in this section we see that there do exist $g\in \SL_2(\mathbb{Z})'$ such that the commutator width is greater than 1.

\begin{defn}
For an admissible value, $t$, \textbf{the minimal width of $t$} is the smallest commutator width of the conjugacy classes with that trace.
\end{defn}

Now how do we find the minimal width of $t$? We start by taking each representatives of the $\SL_2(\mathbb{Z})$ conjugacy classes with trace $t$ (in $\Delta\cap \Gamma(2)'$). If a representative has narrow length 2 or 3 then we know it is a one commutator and can stop. Else, we loop through each representative searching for a mapping that has commutator width 1 and regardless we always take the smallest. If we find one, we are done. We can ignore some relations like transposes of representatives (since if one is a 1-commutator so is the other). To attempt to speed up this search, we sort the matrices (which are already the representatives with the smallest lengths), in increasing order by length, in the hopes that the smallest ones might have a 1-commutator before having to go on to the bigger ones. In the next section we will discuss how to optimize this further.

\textbf{Summary of some search results: } 
\begin{itemize}
  \item First trace that has a class that is not a 1-commutator: 322
  \item First trace that has no classes that are 1-commutators: 322
  \item First trace that has some that are not 1-commutators and some that are: 1170
  \item Traces that have a 1-commutator width representative versus those that have a commutator width 2 representative: up to $|t|<25000$, $267$ had a 1-commutator representative, $160$ had a 2 commutator representative.
  \item First trace that has a 3-commutator: 1298
  \item For $|t|<150000,$ the admissible traces (that aren't failures) have commutator width no more than 2.
\end{itemize}

\begin{conj}\label{1or2comconj}
For admissible values $t$ with $h(t)>0$, the minimal width of $t$ is either 1 or 2.
\end{conj}

\subsection{Connection to the Markoff Equation}

Now, one might ask why we are interested in the commutator width. This interest stemmed from the fact that there is a relationship between 1-commutators and solutions to the level $k$ Markoff-type equations. We will change the wording of some of the theorems stated in Ghosh and Sarnak's paper to fit our purposes. For more information on this topic we recommend \cite{ghosh2017integral}. 

Take the Markoff-type equation for fixed $k,$  $x_1^2+x_2^2+x_3^2-x_1 x_2x_3=k$ and let $M_k$ be the set of integer solutions to this equation. Such things tie together with 1-commutators using Fricke's Identity: 
\begin{equation}\label{eqAB}
  \Tr(A)^2+\Tr(B)^2+\Tr(AB)^2-\Tr(A) \Tr(B)\Tr(AB)=\Tr([A,B])+2.
\end{equation}

We already know that not all of the commutators are 1-commutators (see the summary of the search results) and the algorithm discussed in the previous section by Goldstein and Turner \cite{GoldsteinRichard1979Aotg} appears to find a representative pretty quickly with commutator width 2 (which was the inspiration for Conjecture \ref{1or2comconj}). The goal of this section is to be able to improve the algorithm by using facts about level $k$ Markoff-type solutions to speed up differentiating between 1 and greater than 1 commutator widths.

\begin{defn}
A \textbf{Vieta Involution} of $(x_1, x_2, x_3)$ is defined to be $(x_1, x_2, x_1x_2-x_3).$
\end{defn}

\begin{defn}
An element $x$ is defined to be \textbf{$\hat\Gamma$-equivalent} to $u$ if a combination of permutations, flipping 2 signs and Vieta Involutions on $x$ gives $u$. 
\end{defn}

\begin{lem}[see \cite{ghosh2017integral}]\label{closure}
Let $(a,b,c)$ be an integer solution to the level $k$ Markoff-type equation, $x_1^2+x_2^2+x_3^2-x_1x_2x_3=k,$ and suppose there exists matrices $A,B\in \SL_2(\mathbb{Z})$ such that $\Tr(A)=a,$ $\Tr(B)=b,$ and $\Tr(AB)=c.$ Then any $\hat\Gamma$-equivalent solution, $(x,y,z)$, to $(a,b,c)$ also has matrices in $X,Y\in \SL_2(Z)$ such that $\Tr(X)=x$, $\Tr(Y)=y$, $\Tr(XY)=z$ and $[X,Y]$ is in the same conjugacy class as either $[A,B]$ or $[A,B]^{-1}.$ 
\end{lem}

\begin{proof}
First we take another one of Fricke's identities for $A,B\in \SL_2(\mathbb{R})$: $$\Tr(A)\Tr(B)=\Tr(AB)+\Tr(AB^{-1}).$$ Let $(a,b,c)\in M_k$ and $A,B\in \SL_2(\mathbb{Z})$ such that $\Tr(A)=a$, $\Tr(B)=b$ and $\Tr(AB)=c$. Therefore we can think of $(a,b,c)$ as being mapped to $[A,B].$ We can see that any $\hat\Gamma$- equivalent solution can be mapped to either a matrix in the conjugacy class of $[A,B]$ or its inverse via the following calculations:
\begin{align*}
 (b,a,c)&= (\Tr(B), \Tr(A),\Tr(AB))=(\Tr(B), \Tr(A),\Tr(BA))\mapsto [B,A]=[A,B]^{-1}\\
 (a,c,b)&= (\Tr(A), \Tr(AB),\Tr(B))=(\Tr(A^{-1}),\Tr(AB),\Tr(B))\mapsto [A^{-1},AB]=[A,B]^{-1}\\
 (-a,b,-c)&= (-\Tr(A), \Tr(B),-\Tr(AB))\mapsto [-A,B]=[A,B]\\
 (a,b,a\cdot b -c)&= (\Tr(A), \Tr(B), \Tr(A)\Tr(B)-\Tr(AB))=(\Tr(A), \Tr(B^{-1}),\Tr(AB^{-1}))\\&\mapsto [A,B^{-1}]= B^{-1}[A,B]^{-1}B.
\end{align*}
\end{proof}

The previous lemma is true for all of $\SL_2(\mathbb{Z})$. For $H$ we can invoke the following corollary to see our specific case.
\begin{cor}\label{modingby2}
A solution to the level $k$ Markoff-type equation $x_1^2+x_2^2+x_3^2-x_1x_2x_3=k$ is $\hat\Gamma$-equivalent to a solution $(a,b,c)\equiv (2,2,2)\bmod{4}$ if and only if it itself is congruent to $(2,2,2)\bmod{4}$.
\end{cor}
\begin{cor}\label{closureH}
Let $(a,b,c)$ be an integer solution to the level $k$ Markoff-type equation, $$x_1^2+x_2^2+x_3^2-x_1x_2x_3=k,$$ and suppose there exists matrices $A,B\in H$ such that $\Tr(A)=a,$ $\Tr(B)=b,$ and $\Tr(AB)=c.$ Then any $\hat\Gamma$-equivalent solution, $(x,y,z)$, to $(a,b,c)$ also has matrices in $X,Y\in H$ such that $\Tr(X)=x$, $\Tr(Y)=y$, $\Tr(XY)=z$ and $[X,Y]$ is in the same conjugacy class as either $[A,B]$ or $[A,B]^{-1}.$ 
\end{cor}

\begin{defn}
For our purposes, $k$ is exceptional if $k>0$, $k-2$ is admissible and $k$ is one of the following three forms: (a) $k=u^2+v^2,$ (b) $4(k-1)=u^2+3v^2$, or (c) $k=u^2+4.$ These come from when an element in $M_{k}$ has the absolute value of a term equal to 0, 1, 2 respectively.
\end{defn}

Note that for $k>0$ if we write $k-2=16x+2$ and there exists a solution to the equation $x_1^2+x_2^2+4-2x_1x_2=k,$ then $(x_1-x_2)^2=16x$ which implies that $x$ is a square.  Thus $k-2=16a^2+2$ for some $a.$ So by Lemma \ref{lem2}, we can always write such an admissible value as a 1-commutator with trace $k-2$ and it corresponds to $(2,2,2+4a)\in M_{k}$. So these types of exceptional $k$ we will ignore for the time being.

\begin{thm} \label{biggie} $\text{ }$
\begin{itemize}
  \item Let $t>0$ be admissible and let $t\neq 16a^2+2.$ Set $k=t+2$ and consider the compact set
  \begin{equation*}
    \mathcal{F}_k^+=\{u\in \mathbb{R}^3:\text{ } 6\leq u_1\leq u_2\leq u_3, u_1^2+u_2^2+u_3^2+u_1u_2u_3=k
    \text{ and} u_i \equiv 2 \bmod{4}\}.
  \end{equation*}
  The points in $\mathcal{F}_k^+(\mathbb{Z})=\mathcal{F}_k^+\cap \mathbb{Z}^3$ are $\hat\Gamma$-inequivalent, and any $x\in M_k$ is $\hat\Gamma$-equivalent to a unique point $u'=(-u_1, u_2, u_3)$ where $u=(u_1, u_2, u_3)\in \mathcal{F}_k^+(\mathbb{Z})$ provided that x is of the form $(2,2,2)\bmod{4}.$
  \item Let $t<0$ be admissible. Set $k=t+2$ and consider the compact set 
    \begin{equation*}
    \mathcal{F}_k^-=\{u\in \mathbb{R}^3:\text{ }  6\leq u_1\leq u_2\leq u_3\leq \frac{1}{2}u_1u_2,u_1^2+u_2^2+u_3^2-u_1u_2u_3=k
    \text{ and} u_i \equiv 2 \bmod{4}\}.
  \end{equation*}
   The points in $\mathcal{F}_k^-(\mathbb{Z})=\mathcal{F}_k^-\cap \mathbb{Z}^3$ are $\hat\Gamma$-inequivalent, and any $x\in M_k$ is $\hat\Gamma$-equivalent to a unique point $u=(u_1, u_2, u_3)\in \mathcal{F}_k^-(\mathbb{Z})$ provided that x is of the form $(2,2,2)\bmod{4}.$
\end{itemize}
\end{thm}
\begin{proof}
Suppose that $k$ is not exceptional. Then this is just Theorem 1.1 from \cite{ghosh2017integral} combined with Corollary \ref{modingby2}.

Now suppose that $k$ is exceptional, then by assumption it is of form (a) or (b) from the definition. By Lemma 2.1 of $\cite{ghosh2017integral},$ we know that for any $k$ (not of the form $k=u^2+4$), that any $(x_1,x_2,x_3)\in M_{k}$ falls into one of the following three cases:
\begin{enumerate}
    \item $|x_i|\in \{0,1\}$ for one of the $x_i$'s
    \item $(x_1,x_2,x_3)$ is $\hat\Gamma$-equivalent to an element $(-u_1, u_2, u_3)$ such that $$u_1^2+u_2^2+u_3^2+u_1u_2u_3=k \text{ and } 3\leq u_1\leq u_2\leq u_3$$
    \item $(x_1,x_2,x_3)$ is $\hat\Gamma$-equivalent to $(x_1, x_1x_2-x_3, x_2)$ where $3\leq x_1x_2-x_3\leq x_2\leq x_3$ and $x_1\geq 3.$
\end{enumerate}

By Corollary \ref{modingby2}, we know that any solution with a 0 or a $\pm 1$ in it will not be $\hat\Gamma$-equivalent to a solution congruent to $(2,2,2)\bmod{4}$ and so solutions that fall into case 1 can be ignored. For case 2, we get the statement of the theorem when we restrict to  $(2,2,2)\bmod{4}$ as expected.

Now we look at case 3. In \cite{ghosh2017integral}, a descent argument is used to to prove that for $k\geq 5$ and $|x_i|>2$ we get that any such solution must be $\hat\Gamma$-equivalent to (again) some $(-u_1, u_2, u_3)$ such that $u_1^2+u_2^2+u_3^2+u_1u_2u_3=k$ and $3\leq u_1\leq u_2\leq u_3.$ Thus, we restrict to solutions congruent to  $(2,2,2)\bmod{4}$ and get the statement of the theorem.
\end{proof}

\begin{lem}
Let $g\in \Gamma(2)'$ with $\Tr(g)=t$ and $t\neq 16a^2+2$. Then $g$ is a 1-commutator if only if there exists matrices $A,B\in H$ such that $[A,B]=g$ or $[A,B]^{-1}=g$ and $$(\sgn(-t) \Tr(A),\Tr(B),\Tr(AB))\in\mathcal{F}_{t+2}^{\sgn(t)}(\mathbb{Z}).$$ \end{lem}
\begin{proof}
The one direction of this proof is obvious. Now suppose $g$ is a 1-commutator. Then there exists matrices $X,Y\in H$ such that $[X,Y]=g.$ By Fricke's Identity (\ref{eqAB}), we may conclude that $(\Tr(X), \Tr(Y), \Tr(XY))\in M_{t+2}$ and Theorem \ref{biggie} gives us the existence of some $(a,b,c)$ such that $(\sgn(-t) a,b,c)\in \mathcal{F}_{t+2}^{\sgn(t)}(\mathbb{Z})$ and $(a,b,c)$ is  $\hat\Gamma$-equivalent to $(\Tr(X), \Tr(Y), \Tr(XY)).$  By Corollary \ref{closureH} we know that there then exists $\Tilde{X},\Tilde{Y}\in H$ such that $\Tr(\Tilde{X})=a, \Tr(\Tilde{Y})=b, \Tr(\Tilde{X}\Tilde{Y})=c.$ Furthermore there exists $h\in H$ such that either $[\Tilde{X},\Tilde{Y}]=h[X,Y]h^{-1}$ or $[\Tilde{X},\Tilde{Y}]=h[X,Y]^{-1}h^{-1}.$ Letting $[A,B]=h^{-1}[\Tilde{X},\Tilde{Y}]h,$ the lemma is completed.
\end{proof}

This lemma will help us in two ways. First, when asking what the commutator width of $g\in \Delta \cap \Gamma(2)'$ with $t=\Tr(g)$ and $t$ not of the form $16a^2+2$, we find $\mathcal{F}_{t+2}^{\sgn(t)}(\mathbb{Z}).$ If it is empty, the smallest width it can have is two  and thus we can abort Goldstein and Turner's algorithm if we find a representative that has commutator width 2 given that two is now the best case scenario. Secondly, if $\ell(g)$ is large (and we know that it has a lower bound of $\log(|\Tr(g)|/2)$), Goldstein and Turner's algorithm would take too long and so instead we can try to have Mathematica solve the system of equations that result in an $A,B\in H$ such that for some $(\sgn(-t)a, b, c)\in \mathcal{F}_{t+2}^{\sgn(t)}(\mathbb{Z}),$ $\Tr(A)=a,$ $\Tr(B)=b,$ $\Tr(AB)=c,$ and either $[A,B]=g^{-1}$ or $[A,B]=g.$ If Mathematica found a solution, then $g$ has width 1 and if it proves that there is no solution, then the smallest the commutator width can be is two. However, sometimes Mathematica cannot prove nor disprove the existence and so in that case there are currently no shortcuts. These adaptations to the code are hinted at in the python file from the previous section along with the complementary Mathematica code. \footnote{\href{https://github.com/BrookeOgrodnik/CommutatorSubgroups/blob/master/Mathematica/where\_python\_left\_off.nb}{https://github.com/BrookeOgrodnik/CommutatorSubgroups/blob/master/Mathematica/\newline where\_python\_left\_off.nb}}


\section{Acknowledgements}
The author would like to thank her advisor, Alex Kontorovich, for introducing her to this problem and the reviewer for their extremely helpful feedback. This material is based on work done in the author's PhD dissertation and is partially supported by the National Science Foundation under Grant DMS-1802119.

\bibliographystyle{alpha}
\bibliography{refs.bib}

\end{document}